\documentclass[a4paper,10pt]{article}
\usepackage{amsfonts}
\usepackage{graphicx,amssymb,amsmath}
\usepackage{amscd}
\usepackage[english]{babel}
\usepackage[ansinew]{inputenc}\usepackage{slashbox}
\usepackage[all]{xy}
\usepackage{latexsym}
\usepackage{amsmath}
\usepackage{stmaryrd}

 \hsize \textwidth
 \advance \hsize by -\marginparwidth
\newcommand{\field}[1]{\mathbb{#1}}
\newcommand{\R}{\field{R}}

\newcommand{\N}{\field{N}}

\newtheorem{theo}{Theorem}[section]

\newtheorem{defi}[theo]{Definition}
\newtheorem{lem}[theo]{Lemma}

\newtheorem{prop}[theo]{Proposition}
\newtheorem{rem}[theo]{Remark}

\def\f{\flat}
\def\dis{\ds}
\def\ap{\rightarrow}
\def\S{\Sigma}
\def\l{\lambda}
 \def\a{\alpha}
\def\b{\beta}
\def\g{\gamma}
\def\G{\Gamma}
\def\L{\Lambda}
\def\D{\triangle}
\def\t{\tau}
\def\d{\delta}
\def\o{\omega}

\def\l{\lambda}
\def\n{\nu}
\def\f{\flat}

\def\r{\rho}

\def\n{\nu}
\def\s{\sigma}

\def\so{\underline}

\def\O{\Omega}

\def\e{\epsilon}

\def\~{\tilde}
\def\dis{\displaystyle}
\input epsf

\title{{ Snakes and  articulated arms in an  Hilbert space}}
 \author{ F. Pelletier \&  R. Saffidine}

\date{}

\begin{document}

\maketitle

\begin{abstract}
The purpose of this paper is to give an illustration of results on integrability of distributions  and orbits of vector fields on Banach manifolds obtained in \cite{Pe} and \cite{LaPe}.
Using  arguments and results of these  papers, in the context of a separable Hilbert  space,
 we give a generalization of a  Theorem of accessibility   contained in \cite{Ha}, \cite{Ro} and proved  for a finite dimensional Hilbert space.
    \end{abstract}

\section{Introduction}
In finite dimension,  a  snake (of length L) is a (continuous) piecewise $C^1$-curve $S : [0,L] \ap  \R^d$, arc-length parameterized  such that S(0) = 0. According to \cite{Ha} for articulated arms  (i.e. when $S$ is affine by parts) and \cite {Ro} in the general case, "charming  a snake" is a control problem so that   its  "head"   $S(L)$ describes   a given $C^1$-curve $c: [0,1]\ap \R^d$ in minimal way. More precisely
we look for a $1$-parameter family $\{S_t\}_{t\in [0,1]}$ such that $S_t(L)=c(t)$ for all $t\in [0,1]$ so that the family $\{S_t\}$ has a minimal infinitesimal kinematic energy. We can formulate this problem in the following way:

Each snake $S$ of length $L$ in $\R^d$ can be given by a  piecewise $C^0$-curve $u : [0,L] \ap  \mathbb{S}^{d-1}$ such that $S(t)=\int_0^tu(\t)d\t$. We look for a $1$-parameter  family $\{u_t\}_{t\in [0,1]}$  such that the associated  family $S_t$ of snakes satisfies  $S_t(L)=c(t)$ for all $t\in [0,1]$ so that the infinitesimal kinematic energy $\dis\frac{1}{2}\int_0^L||\frac{d}{dt}u_t(s)||ds$ is minimal.\\

The purpose of this paper is to give a generalization of this problem in the context of separable Hilbert spaces. More precisely, given a separable Hilbert space $\field{H}$ we consider the smooth hypersurface  $\field{S}^\infty$ of elements of norm $1$. As previously, a Hilbert snake of length $L$ is a continuous piecewise
$C^1$-curve $S : [0,L] \ap  \field{H}$, arc-length parameterized  such that  $S(0) = 0$. An articulated arm corresponds to the particular case where $u$ is affine in each part.  Then a snake is also given  by a  piecewise $C^0$-curve $u : [0,L] \ap  \field{S}^\infty$ such that $S(t)=\int_0^tu(\t)d\t$.  Given a fixed partition $\cal P$ of $[0,L]$, the set  ${\cal C}^L_{\cal P}$ of such curves will be called the configuration set and carries a natural structure of Banach manifold. For articulated arms, the configuration space is the subset ${\cal A}^L_{\cal P}$ of $u$ which are constant on each subinterval associated to the partition.  In fact, ${\cal A}^L_{\cal P}$ is a weak Hilbert submanifold of ${\cal C}^L_{\cal P}$.\\

To any "configuration" $u\in {\cal C}^L_{\cal P}$ is naturally associated the "end map"
${\cal E}(u)=\dis\int_0^Lu(s)ds$. This map is smooth and its kernel has a canonical complemented subspace which gives rise to a closed distribution $\cal D$ on ${\cal C}^L_{\cal P}$. In finite dimension,  for a one parameter  family $\{u_t\}_{t\in [0,1]}$   the associated  family $S_t$ of snakes satisfies  $S_t(L)=c(t)$ for all $t\in [0,1]$ so that the infinitesimal kinematic energy $\dis\frac{1}{2}\int_0^L||\frac{d}{dt}u_t(s)||ds$ is minimal. If $c(t)$ has a "lift "  $\tilde{c}$ in ${\cal C}^L_{\cal P}$ which is tangent to $\cal D$, it is called a "horizontal lift". So  the problem for the head of the Hilbert snake  to join an initial state $x_0$  to a final state $x_1$ can be transformed in the following "accessibility problem":

Given an initial (resp. final) configuration $u_0$ (resp. $u_1$) in ${\cal C}^L_{\cal P}$, such that ${\cal E}(u_i)=x_i$, $i=0,1$, find a piecewise $C^1$ horizontal curve $\g:[0,T]\ap {\cal C}^L_{\cal P}$ (i.e. $\g$ is tangent to $\cal D$) and  which joins $u_0$ to $u_1$. \\

So, given any configuration $u\in {\cal C}^L_{\cal P}$ we look for the accessibility set ${\cal A}(u)$ of  all configurations $v\in {\cal C}^L_{\cal P}$ which can be joined from $u$ by a piecewise $C^1$ horizontal curve. In the context of finite dimension, in \cite{Ha} and  \cite{Ro},  using arguments about  the   action of the Mo\"{e}bus group on ${\cal C}^L_{\cal P}$, it can be shown that  ${\cal A}(u)$  is the maximal integral manifold of a finite dimensional distribution on ${\cal A}^L_{\cal P}$  and  ${\cal C}^L_{\cal P}$. Unfortunately, in our context, the same argument does not work. Moreover, as we are in the context of infinite dimension for $\field{S}^\infty$, we cannot hope to get a finite dimensional distribution whose maximal integral manifold is ${\cal A}(u)$.\\

However, our principal result is to construct a canonical distribution $\bar{\cal D}$ modeled on a Hilbert space, which is integrable and such that the accessibility set ${\cal A}(u)$ is a dense subset  of  the maximal integral manifold through $u$ of $\bar{\cal D}$. Moreover this distribution is minimal in some natural sense (see Remark \ref{min}). In fact, when $\field{H}$ is finite dimensional, $\bar{\cal D}$ is exactly the finite distribution obtained in \cite{Ro} whose leaves are the accessibility sets.\\

The arguments used in our proof can be found in \cite{Pe} and \cite{LaPe}. Moreover, this Theorem of accessibility can be seen as an application of results obtained  in \cite{LaPe}; it is also an illustration of the almost Banach algebroid structures developed in \cite{CaPe} (see subsection \ref{ALalgebroid}).\\

This paper is organized as follows.  Section \ref{rappel} contains all definitions and results of \cite{Pe} and \cite{LaPe} which are used in the proof about  the accessibility sets. {\it In a first time, the reader can skip this  section; he can  only  refer to this paragraph for a deeper reading}. In  section 3, we define the configuration space, its Banach manifold structure and we construct the horizontal distribution. The
last section  presents in more detail the previous optimal problem and contains the principal result (Theorem \ref{acesset} in subsection \ref{result}). The proof of this Theorem which needs all definitions and results recalled in section \ref{rappel} is developed in subsections \ref{barD} and \ref{access}.
\section{Preliminaries}\label{rappel}
\subsection{Weak distributions on a Banach manifold}\label{prelim}
 In this subsection, from\cite{Pe} we recall all definitions, properties and results we shall use later.\\

 Let  $M$ be a connected Banach manifold modeled on a Banach space $E$. We denote by ${\cal X}(M)$ the set of local vector fields on $M$. The flow of  any $X\in {\cal X}(M)$ will be denoted by $\phi^X_t$. We then have the following definitions and properties: \\

$\bullet$ A {\bf weak submanifold}  of $M$  is a pair $(N,f)$ where $N$  is a connected  Banach manifold   (modeled on a Banach space $F$) and $f:N\ap M$ is a smooth  map such that : 

--- there exists a continuous injective linear  map $i:F\ap E$ between these two Banach spaces;

--- $f $ is injective and the tangent map $T_xf:T_xN\ap T_{f(x)}M$ is injective for all $x\in N$.\\

 Note that for  a weak submanifold $f:N\ap M$, on the subset $f(N)$ of $M$ we have two topologies:

---  the induced topology from $M$;

---  the topology for which $f$ is a homeomorphism from $N$ to $f(N)$.\\
With this last topology, via $f$, we get a structure of Banach manifold modeled on $F$.
 Moreover, the inclusion from $f(N)$ into $M$ is continuous as a map from the Banach manifold $f(N)$ to $M$.
 In particular, if $U$ is an open set of $M$, then $f(N)\cap U$ is an open set for the topology of the Banach manifold on $f(N)$.\\

$\bullet$ According to \cite{Pe}, a  {\bf weak distribution} on $M$  is an assignment  $ {\cal D}: x \mapsto {\cal D}_x$  which, to every $x\in M$, associates  a vector subspace  ${\cal D}_x$ in $T_xM$  (not necessarily closed) endowed with a norm $||\;||_x$ such that   $({\cal D}_x, ||\;||_x)$ is a Banach space (denoted by $\tilde{\cal D}_x$) and such  that the natural inclusion $i_x : \tilde{\cal D}_x \ap T_xM$ is continuous. Moreover, if the Banach structure on ${\cal D}_x$ is a Hilbert structure, we say that $\cal D$ is a {\bf weak Hilbert distribution}.\

When ${\cal D}_x$ is closed, we have a natural Banach structure  on $\tilde{\cal D}_x$, induced by the Banach structure on $T_xM$, and so we get the classical definition of a distribution;  in this case we will say that  $\cal D$ is {\bf closed}.

\noindent A  (local) vector field $Z$ on $M$ is {\bf tangent} to $\cal D$, if for  all $x\in$ Dom($Z$), $Z(x)$ belongs to ${\cal D}_x$. The set of local vector fields tangent to $\cal D$  will be denoted by {\bf $\bf {\cal X}_{\cal D}$}.

$\bullet$ We say that  $\cal D$ is {\bf  generated by a subset  }  ${\cal X}\subset {\cal X}(M)$ if, for every $x\in M$, the vector space ${\cal D}_x$ is the  linear hull  of the set $\{Y(x)\;,\;Y\in{\cal X}\;,\;x\in$ Dom$(Y)\}$.  \\

For a weak distribution $\cal D$ on $M$ we have the following definitions:\\

$\bullet$ an {\bf integral manifold} of $\cal D$  through $x$ is a weak  submanifold $f:N\ap M$ such that there exists $u_0\in N$ such that $f(u_0)=x$ and $T_uf(T_uN)={\cal D}_{f(u)}$ for all $u\in N$.\\

$\bullet$   $\cal D$ is called {\bf integrable} if for any $x\in M$ there exists an integral manifold $N$ of $\cal D$ through $x$.\\

$\bullet$ if  $\cal D$ is generated by a set $\cal X$ of local vector fields, then  $\cal D$ is  called {\bf ${\bf {\cal X}}$- invariant} if for any   $X\in {\cal X}$, the tangent map $T_x\phi^X_t $ sends ${\cal D}_x$ onto $ {\cal D}_{\phi^X_t(x)}$ for all $(x,t)\in \O_X$. ${\cal D}$ is {\bf invariant} if ${\cal D}$ is $ {\cal X}_{\cal D}-$ invariant. \\

{\it Now we introduce essential properties of "local triviality" which will play an essential role thorough this paper}:\\

$\bullet$   $\cal D$ is {\bf lower (locally) trivial} if for each $x\in M$, there exists an open neighborhood $V$ of $x$,  a smooth map $\Theta:\tilde{\cal D}_x\times  V \ap TM$  (called {\bf  lower trivialization}) such that :
\begin{enumerate}
\item[(i)]  $ \Theta(\tilde{\cal D}_x\times\{y\})\subset {\cal D}_y$ for each $y\in V$

\item[(ii)] for each $y\in V$,  $\Theta_y\equiv \Theta(\;,y):\tilde{\cal D}_x\ap T_yM$ is a continuous operator  and $\Theta_x:\tilde{\cal D}_x\ap T_xM$  is the natural inclusion $i_x$

\item [(iii)] there exists a  continuous operator $\tilde{\Theta}_y: \tilde{\cal D}_x\ap \tilde{\cal D}_y$ such that $i_y\circ \tilde{\Theta}_y=\Theta_y$, $\tilde{\Theta}_y$ is an isomorphism from $\tilde{\cal D}_x$ onto ${\Theta}_y(\tilde{\cal D}_x)$
and  $\tilde{\Theta}_x$ is the identity of $\tilde{\cal D}_x$\\
\end{enumerate}

$\bullet$  ${\cal D}$ is called  {\bf (locally) upper  trivial} if,  for each $x\in M$, there exists an open neighborhood $V$ of $x$, a Banach space ${F}$ and  a smooth map $\Phi:F\times  V \ap TM$   such that :
\begin{enumerate}
\item[(i)]  $ \Phi(F\times\{y\})= {\cal D}_y$ for each $y\in V$
\item[(ii)] for each $y\in V$,  $\Phi_y\equiv \Phi(\;,y):F\ap T_yM$ is a continuous operator such that $\Phi_y(F)={\cal D}_y$
\end{enumerate}

$\bullet$   ${\cal D}$  is called {\bf strong upper  trivial} if, for any $x\in M$, there exists an upper trivialization $\Psi:F\times  V \ap TM$ such
 $\ker \Psi_x$ is complemented (i.e. $F=\ker \Psi_x \oplus S$) such that the restriction $\theta_y$  of $\Psi_y$ to $S$  is injective for any $y\in v$,  and then ${\Theta}(u,y)=({\theta}_y\circ [{\theta}_x]^{-1}(u), y)$ is  a lower trivialization of ${\cal D}$. In this case $\Theta$ is called the {\bf associated lower trivialization}.\\

A strong  upper trivial weak distribution ${\cal D}$ is called {\bf Lie bracket invariant}  if, for any $x\in M$, there exists an upper trivialization $\Phi:F\times V\ap TM$ such that for any $u\in F$ , there exists $\varepsilon>0$,  such that, for all $0<\t<\varepsilon$,
 we have a smooth field of operators $C:[-\t,\t]\ap L(F,F)$ with the following property
 \begin{eqnarray}\label{condLieinvst}
 [X_u,Z_v](\g(t))=\Phi(C(t)[v],\g(t)) \textrm{ for any } Z_v=\Phi(v,\;) \textrm{ and any  } v\in F
\end{eqnarray}
along the integral curve $\gamma : t\mapsto \phi^{X_u}_t(x)$ on $[-\t,\t]$ of  the lower section $X_u=\Theta(\Phi(u,x),\;)$ .

With these definitions we have:

\begin{theo}\label{lieinv}${}$\\
Let ${\cal D}$ be a strong upper trivial weak  distribution. Then ${\cal D}$ is integrable if and only if ${\cal D}$ is Lie bracket invariant.\\
\end{theo}

\subsection{Orbit of a family of vector fields}\label{orbit}
In this subsection we expose the results of \cite{LaPe} which will be useful for the proof of Theorem \ref{acesset}.\\

Let  ${\cal X}$ be a set of local vector fields on $M$.
 Given $x\in M$, we say that
 $\cal X$ satisfies the condition
(LB(s))  at $x$ (Locally Bounded of order $s$), if there exists a chart
$(V_x, \phi)$ centered at $x$ and a constant $k>0$ such that:\\
for any $X\in {\cal X} $,  whose  domain
Dom$(X)$ contains $V_x$, we have
\begin{eqnarray}\label{LBs}
\sup \{||J^s[\phi_*X](y)||,\; X\in {\cal X} , \; y\in
V_x\}\leq k.
\end{eqnarray}
For any finite or  countable  ordered set  $A$ of indexes, consider   a family $\xi=\{X_\a\}_{\a\in A}$ where the $X_\a$ are defined on a same open set $V$ and satisfies the condition (LBs) for $s\geq 1$. Given any bounded integrable map $u=(u_\a)_{\a\in A} $ from some interval $I$ to $l^1(A)=\{\t=(\t_\a),\; \dis\sum_{\a\in A}|\t_\a|<\infty\}$ we can associate
a time depending  vector field of type
$$Z(x,t,u)=\dis\sum_{\a \in A}u_\a(t)X_\a(x),$$
For such a vector field
there exists a {\bf flow $\Phi^\xi_u(t,\;)$} (see Theorem 2 of \cite{LaPe}). \\
Given some $\t\in l^1(A)$,  we set $||\t||_1=\dis\sum_{\a\in A}|\t_\a|$. On the corresponding interval $[0,||\t||_1$], we consider the partition $(t_\a)_{\a\in A}$ of  this interval defined by, $t_0=0$ and  for $\a\in A$, $t_\a=\dis\sum_{\b=1}^\a |\t_\b|$.
If we choose   $u=\G^\t=(\G^\t_\a)$ where $\G^\t_\a$ is the indicatrix function of  $]t_{\a},t_{\a+1}[$ we can associate to $(\xi,\t )$ a time depending vector field $Z(x,t,u)$ as previously. Under appropriate assumptions, for such a $Z$,  we get an associated  flow, denoted by $\Phi^\xi_\t(t,\;)$ .
Assume that the set of  all Dom$X^($\footnote{Dom$X$ is  the maximal open set on which $X$ is defined}$^)$  for $X\in{\cal X}$ is a covering of $M$ and is bounded at each point, i.e.   the set of values $\{X(x),\; X\in {\cal X}\}\subset T_xM$ is bounded for any $x\in M$. We can  enlarge $\cal X$  to the set  $\hat{\cal X}$ given  by
$$\hat{\cal X}=\{Z=\Phi_*(\n Y),\; Y\in {\cal X},\; \Phi=\phi^{X_p}_{t_p}\circ\cdots\circ\phi^{X_1}_{t_1} \textrm{ for } X_1,\cdots, X_p\in {\cal X} ; \textrm{  and appropriate }\n\in \R\} $$
(see subsection 3.1 of \cite{LaPe}). Then  $\hat{\cal X}$ satisfies the same previous properties as $\cal X$. From this set $\hat{\cal X}$, we associate  an appropriate pseudo-group ${\cal G}_{\cal X}$ of local diffeomorphisms which are finite compositions of  flows of type $\phi^X_t$ with $X\in{\cal X}$ and   of type $\Phi_u^\xi(||\t||_1,.)$ (as we have seen previously) or its inverse for $\xi \subset\hat{\cal X}$.

\noindent  To ${\cal G}_{\cal X}$ is naturally associated the following equivalence relation
on  $M$:
$$x\equiv y \textrm{ if and only if there exists } \Phi\in{\cal G}_{\cal X} \textrm{ such that } \Phi(x)=y$$

\noindent {\bf An equivalence class is called a $\bf {\cal X}$-orbit}.\\

\begin{prop}\label{l1Xorbit}\cite{LaPe}${}$
 For each pair $(x,y)$ in the same $ {\cal X}$-orbit  either we have a piecewise smooth curve  which joins $x$ to $y$ and whose each smooth part is tangent to $X$ or $-X$ for some $X\in {\cal X}$  or there exists a sequence $\g_k$ of such  piecewise smooth curves whose origin is $x$ (for all curves) and whose sequence of ends converges to $y$.
\end{prop}
 Consider any set $\cal Y$ of local vector fields which contains $\hat{\cal X}$. Assume that there exists a weak distribution $\D$ generated by $\cal Y$ which  is integrable on $M$ and for each $x\in M$
there exists a lower trivialization $\Theta :F\times V\ap TM$  for some Banach space $F$ (which depends of $x$) and some neighborhood $V$ of $x$ in $M$.  Let $N$ be the union of all integral manifolds $i_L:L\ap M$ through $x_0$. Then $i_N:N\ap M$ is the maximal integral manifold of $\D$ through $x_0$(see Lemma 2.14 \cite{Pe}).

  \begin{prop}\label{varXS}${}$ (see \cite{LaPe})\\
  As previously, let $f: N\ap M$ be the maximal integral manifold of $\D$  through $x$.
  \begin{enumerate}
 \item Let $Z\in {\cal X}(M)$ be such that Dom$(Z)\cap  f(N)\not=\emptyset$  and  $Z$ is tangent to $\D$. Set $\tilde{V}_Z=f^{-1}(\textrm{Dom}(Z)\cap  f(N))$. Then $\tilde{V}_Z$ is an open set in $N$ and there exists a vector field $\tilde{Z}$ on  $N$ such that Dom$(\tilde{Z})=\tilde{V}_Z$ and $ f_*\tilde{Z}=Z\circ f $.

Moreover, if $]a_x,b_x[$ is the maximal interval on which the integral curve $\g: t \mapsto \phi^Z(t,{x})$ is defined  in $M$, then the integral curve $\tilde{\g}:t\ap \phi^{\tilde{Z}}(t,\tilde{x})$ is also defined on  $]a_x,b_x[$ and we have
\begin{eqnarray}\label{gZ}
\g=f\circ\tilde{\g}
\end{eqnarray}
\item Let be $\xi=\{X_\b,\; \b\in B\}\subset \hat{\cal X}\subset {\cal Y}$ which satisfies the conditions (LB(s)) on a chart domain $V$ centered at $x\in f(N)$ and consider  the associated flow $\Phi^\xi_\t$.
For some $\t\in l^1(B)$ let  $\g$ be the curve on $[0,||\t||_1]$ defined by $\g(t)=\Phi^\xi_\t(t,x)$. Then  there exists a curve $\tilde{\g}:[0,||\t||_1[\ap N$ such that
\begin{eqnarray}\label{l1N}
f\circ \tilde{\g}=\g \textrm{ on } [0,||\t||_1[
\end{eqnarray}
\end{enumerate}
\end{prop}

According to the properties of $\cal X$  we can associate to this set a weak distribution $\cal D$  in the following way:

${\cal D}_x=\{Y=\dis\sum_{X\in{\cal X}}\l_XX(x)\}$ for any  absolutely summable family $\{\l_X,\; X\in {\cal X},\; x \in \textrm{Dom}(X)\}$

In the same way we can also associate to $\hat{\cal X}$ a weak distribution $\hat{\cal D}$ which contains $\cal X$ and which is ${\cal X}$-invariant. Moreover,  for a set $\cal Y$ of local vector fields which contains ${\cal X}$ and which is bounded at each point, we can also associate a weak distribution $\D$ of the previous type.
If $\D$ is $\cal X$ invariant, then $\hat{\cal D}_x\subset \D_x$ for any $x\in M$.

To the set ${\cal X}$ we can associate the sequences of families
 $${\cal X}={\cal X}^1\subset {\cal X}^2={\cal X}\cup\{[X,Y],\; X,Y\in {\cal X}\}\subset \cdots\subset {\cal X}^k={\cal X}^{k-1}\cup\{[X,Y],\; X\in {\cal X}, Y\in {\cal X}^{k-1}\}\subset\cdots$$
When ${\cal  X}^k$ is bounded at each point, as previously, we can associate a weak distribution ${\cal D}^k$ generated by ${\cal X}^k$.\\

Consider an ordered finite or countable set of indexes $A$ and assume that we have  $\hat{\cal D}$ fulfilling the following conditions
\begin{enumerate}
\item  for any $x\in M$ there exists a strong upper trivialization $\Phi:l^1(A)\times V\ap TM$ such that $\Phi(e_\a,.)=Y_\a(.)$ for each $\a\in A$ where $\{e_\a\}_{\l\in A}$ is the canonical basis of $l^1(A)$;
\item for any $x\in M$ there exists   a neighborhood $V$ of $x$  such that, $V\subset \dis\cap_{\a\in  A}\textrm{Dom}(Y_\a)$,  and a constant $C>0$ such that  we have
\begin{equation}\label{hypinvol}
[Y_\a,Y_\b](y)=\dis\sum_{\n\in A}C_{\a\b}^\n(y) Y_\n(y) \textrm{ for any } \a,\b\in A
\end{equation}
where each $C_{\a\b}^\n$ is a smooth function  on $V$, for any  ${\a,\b,\n\in A}$ and  we have $$\dis\sum_{\a,\b,\n\in A}|C_{\a\b}^\n(y)|\leq C$$ for any $y\in V$.
 \end{enumerate}

 Then we have:

 \begin{theo}\label{III"} (\cite{LaPe}):
 \begin{enumerate}
 \item Under the previous assumptions,  the distribution $\hat{\cal D}$  is integrable and
each $\cal X$-orbit $\cal O$ is the union of the maximal integral manifolds which meet  $\cal O$ and such an integral manifold is dense in $\cal O$.
  \item  If ${\cal D}^k$  is defined  and satisfies  the previous assumptions   for some  $k\geq 2$, then  we have ${\cal D}^k=\hat{\cal D}$  and ${\cal D}^k$ is integrable.
\end{enumerate}
 \end{theo}





\section{ Hilbert snakes and Hilbert  articulated  arms}
\subsection{The configuration space}
Let  $\mathbb{H}$ be a separable Hilbert space  and $<.,.>$ (resp. $||.||$) the inner product (resp. the norm) on $\mathbb{H}$.
 We consider a fixed hilbertian basis $\{ e_{i}\}_{i\in \mathbb{N}} $ in  $\mathbb{H}$. Any  $ x\in\mathbb{H}$ will be written as a serie $x=\dis \sum_{i\in\mathbb{N}}x_{i} e_{i}$
where $x_{i}=\langle x,e_{i} \rangle$ is the $i^{th}$coordinate of $x$. We denote by
$\mathbb{S}^{\infty }=\left\{ x\in\mathbb{H}\;:\;\left\Vert x\right\Vert =1\right\} $ the unit sphere in $\mathbb{H}$. Note that $\mathbb{S}^{\infty }$ is a codimension one hypersurface in $\mathbb{H}$ whose smooth equation is $||x||^2-1=0$.\\
A curve $\gamma :\left[ a,b\right] \rightarrow M$ called ${C}
^{k\text{ }}$ piecewise if there exists  a finite set $\mathcal{P}$=$
\left\{ a=s_{0}<s_{1}<...<s_{N}=b\right\} $ such that, for all $i=0,...,N-1$,
the restriction of  $\gamma$ to the interval $\left[ s_{i},s_{i+1}\right[$ can be extended to a curve  of class ${C}^{k}$ on the closed interval
$\left[ s_{i},s_{i+1}\right] $.

Given any metric space $(X,d)$, we denote by  $\mathcal{C}
\left( \left[ a,b\right] ,X\right) $  the set of  continuous curves $u:[a,b]\ap X$. Recall that on  $\mathcal{C}
\left( \left[ a,b\right] ,X\right) $ we have the usual distance  $d_\infty$ defined by
$$d_\infty(u_1,u_2)=\dis\sup_{t\in [a,b]} d(u_1(t),u_2(t))$$
and ${\cal C}\left( \left[ a,b\right] ,X\right),d_\infty) $ is a complete metric space.\\

 For a given partition $\mathcal{P}$=$
\left\{ a=s_{0}<s_{1}<...<s_{N}=b\right\} $ of $[a,b]$,
 let be $\mathcal{C}_{\mathcal{P}}^k\left( \left[
a,b\right] ,\mathbb{S}^\infty\right) $ (resp. $\mathcal{C}_{\mathcal{P}}^k\left( \left[
a,b\right] ,\mathbb{H}\right) $ the set of curves $u\in\mathcal{C}
\left( \left[ a,b\right] ,\mathbb{S}^\infty\right)$ (resp. $u\in\mathcal{C}
\left( \left[ a,b\right] ,\mathbb{H}\right)$)
 which are $C^k$-piecewise relatively to $\mathcal{P}$ for $k\in \N$.\\

\noindent Thorough this paper, we fix a real number $L>0$  and  $\mathcal{P}$ is a given fixed partition of $[0,L]$. \\

{\bf A Hilbert snake} is a continuous  piecewise $C^1$-curve $S : [0,L] \ap  \field{H}$, such that $||\dot{S}(t)||=1$ and S(0) = 0. When $S$ is affine by part, we  call this snake {\bf an affine snake or a Hilbert articulated arm}.

In fact, a snake is characterized by $u(t)=\dot{S}(t)$ and of course we have $S(t)=\dis\int_0^tu(s)ds$  where $u:[0,L] \ap  \field{S}^\infty$ is a piecewise $C^0$-curve associated to the partition $\cal P$. Moreover, this snake is affine  if and only if $u$ is constant on each subinterval of $\cal P$.\\

The set
$${\cal C}^L_\mathcal{P}={\mathcal{C}}_{\mathcal{P}}^0\left(
\left[ 0,L\right] ,
\mathbb{S}^\infty \right) .$$

is called the {\it configuration space of the snakes}  in $\mathbb{H}$ of length $L$ relative to the partition $\cal P$ . We can also put on ${\cal C}^L_\mathcal{P}$ the distance $d_\infty$ defined by
$$d_\infty(u_1,u_2)=\dis\sup_{t\in [a,b]} ||(u_1(t),u_2(t))||$$

Note that the subset
$${\cal A}^L_{\cal P}=\{u\in {\cal C}^L_{\cal P}, \textrm{ such that } u \textrm{ is constant on each subinterval } [s_{i-1},s_i[, i=1,\cdots N\}$$
 is the configuration space of Hilbert articulated arms  in $\mathbb{H}$ of length $L$ relative to the partition $\cal P$.

The natural map
\begin{eqnarray}\label{homeo}
  h: && {\cal C}^L_{\cal P}\rightarrow\prod_{i=0}^{N-1}{\cal C}^{0}([s_{i},s_{i+1}],\field{S}^{\infty})\nonumber\\
   && u\mapsto(u\mid_{[s_{0},s_{1}]},...,u\mid_{[s_{i},s_{i+1}]},u\mid_{[s_{N-1},s_{N}]})
\end{eqnarray}
is a homeomorphism. In particular,  $({\cal C}^L_\mathcal{P},d_\infty)$ is a complete metric space. Note that the restriction of $h$ to ${\cal A}^L_{\cal P}$ is a homeomorphism onto $[\field{S}^\infty]^N$. Moreover we have

\begin{prop}\label{banachst} ${}$\\  ${\cal C}^L_\mathcal{P}$ has a structure of Banach manifold and according to (\ref{homeo}) the natural  map  $$h:  {\cal C}^L_{\cal P}\rightarrow\prod_{i=0}^{N-1}{\cal C}^{0}([s_{i},s_{i+1}],\field{S}^{\infty})$$ is a diffeomorphism. Moreover
${\cal A}^L_{\cal P}$ is a weak Hilbert  submanifold diffeomorphic to $[\field{S}^\infty]^N$ and the topology associated to this structure  and the topology induced by  ${\cal C}^L_{\cal P}$ coincide.
\end{prop}

\textit{Proof }:  as ${\cal C}^L_{\cal P}$ is homeomorphic to $\prod_{i=0}^{N-1}C^{0}([s_{i},s_{i+1}],\field{S}^{\infty})$ it is sufficient to prove that ${\cal C} ^L_{\cal P}$ has a natural structure of Banach manifold ${\cal C}([0,L], \field{S}^\infty)$for ${\cal P}=\{0,L\}=$. Note that ${\cal C}([0,L], \field{H})$ is a Banach space for the norm
$$||u||_\infty=\dis\sup_{t\in [0,L]}||u(t)||$$
On the other hand,   consider the map $\eta: {\cal C}([0,L], \field{H})\ap {\cal C}([0,L],\R)$ defined by $\eta(u)(t))=||u(t)||^2$. Of course we have ${\cal C}([0,L], \field{S}^\infty)=\eta^{-1}(1)$.  As the function $x\ap ||x||^2$ is analytic on $\field{H}\setminus\{0\}$, the map  $\eta$  is also analytic on a neighborhood of  ${\cal C}([0,L], \field{S}^\infty)$ and its differential is \\
$D_{u}\eta(v)(s)=2\langle u(s),v(s)\rangle$. To end the proof, it is sufficient to prove that  $D\eta$ is surjective on ${\cal C}([0,L], \field{S}^\infty)$ and at each $u\in {\cal C}([0,L], \field{S}^\infty)$ its kernel is complemented. As in finite dimension (see \cite{Ro}) , for  $u\in {\cal C}([0,L], \field{S}^\infty)$ we have
\begin{enumerate}
\item  for any $f\in {\cal C}([0,L],\R)$,  for $u\in {\cal C}([0,L], \field{S}^\infty)$ we have $D_u\eta(v)=f$ for
$$v(s)=\dis\frac{1}{2}f(s)\dis\sum_{i\in \N} u_i(s)e_i$$
\item
 $ \ker D_{u}\eta= \{ v \in{ \cal C}([0,L],\mathbb{H}) \mid \langle u(s),v(s)\rangle=0 \forall s\in[0,L]\}$

Consider the closed subspace
    $H_{u}=\{f(s)u(s) \textrm{ for } f\in{\cal  C}([0,L],\mathbb{R})\}$
of ${\cal C}([0,L], \field{H})$.
Then each  $v \in {\cal C}([0,L],\mathbb{H}) $ can be written as
 $$ v(s)=\underbrace{v(s)- \langle u(s),v(s)\rangle u(s)}_{\in \ker D_{u}\eta } +\underbrace{ \langle u(s),v(s)\rangle u(s)}_{\in H_{u}}$$
Of course clearly $H_{u}\cap \ker D_{u}\eta={0}$
\end{enumerate}

On the other hand,  when ${\cal P}=\{0,L\}$, the map $h$ is nothing but the identity, so, according to the product structure, we see that $h$ is a smooth diffeomorphism.

Again, when ${\cal P}=\{0,L\}$,  the map $f:\field{S}^\infty\ap  {\cal A}^L_{\cal P}$ defined by $x\ap f(x)(t)=x$ is a homeomorphism from $\field{S}^\infty$ to ${\cal C}^L_{\cal P}$ which is the restriction of the map  $\bar{f}:\field{H}\ap {\cal C}^0([0,L],\field{H})$ defined in the same way. But $\bar{f}$ is linear and injective  and then smooth and its differential is also injective. It follows that the restriction of $\bar{f}$ to $\field{S}^\infty$ has the same properties  as the image of this restriction is precisely ${\cal A}^L_{\cal P}$; we conclude that ${\cal A}^L_{\cal P}$ is a weak submanifold of ${\cal C}^L_{\cal P}$ diffeomorphic to $\field{S}^\infty$.  Moreover, as the canonical  topology on  $\field{H}$ coincides with the topology induced by $||\;||_\infty$, the same is true for each  induced  topology on ${\cal A}^L_{\cal P}$. According to the product structure of manifolds, this ends the proof of Lemma \ref{banachst}.

${}$\hfill $\D$\\

The tangent space $T_u{\cal C}^L_{\cal P}$
    can be identified with the set
$$\{ v\in {\cal C}^0_{\cal P}([0,L],\field{H}) \textrm { such that } <u(s),v(s)>=0 \textrm{ for all } s\in [0,L]\}$$
 This space is naturally provided with the induced  norm $||.||_\infty$.  On the other hand, note  that any $v\in {\cal C}^0_{\cal P}([0,L],\field{H})$ is integrable on $[0,L]$ and so we get a inner product on this space given by:

 \begin{eqnarray} \label{prod}
<v,w>_{L^2}=\dis\int_{0}^L<v(s),w(s)>ds
\end{eqnarray}

This inner product induces a natural norm $||.||_{L^2}$  on $T_u{\cal C}^L_{\cal P}$  given by:
$$||v||_{L^2}=[\dis\int_{0}^L<v(s),v(s)>ds]^{\frac{1}{2}}$$
We have the following inequality
\begin{eqnarray}\label{inegalitenorm}
||u||_{L^2}\leq \sqrt{L}||u||_\infty.
\end{eqnarray}

In the same way the tangent space  $T_u{\cal A}^L_{\cal P}$ can be  identified with the set $v=(v_1,\cdots,v_N)\in \field{H}^N$ such that $<v_i,u_i>=0$ for $i=1\cdots N$ if $u=(u_1,\cdots,u_N)$. Of course, this vector space can be also considered as a subspace of $T_u{\cal C}^L_{\cal P}$. Note that this subspace in closed in $T_u{\cal C}^L_{\cal P}$.

\begin{rem}\label{weakmetric}${}$\\
As $(T_u{\cal C}^L_{\cal P}, ||.||_{L^2})$ is not complete, the  norm $||.||_\infty$ and $||.||_{L^2}$ are not equivalent (on each $T_u{\cal C}^L_{\cal P}$). So the inner product defined by (\ref{prod}) gives rise only to a weak Riemannian  ${G}$ metric on $T{\cal C}^L_{\cal P}$.\\
 As ${\cal A}^L_{\cal P}$ is diffeomorphic to $[\field{S}^\infty]^N$, the tangent  $T_u{\cal A}^L_{\cal P}$  can be identified with
 $$T_{x_1}\field{S}^\infty\times\cdots\times T_{x_N}\field{S}^\infty$$
 for $u=(x_1,\cdots,x_N)\in [\field{S}^\infty]^N$. So the canonical inner product
  on $\field{H}^N$, induces an natural inner product on  $T_u{\cal A}^L_{\cal P}$.\\
    On the other hand,  the inner product $<\;,\;>_{L^2}$ on $T_u{\cal C}^L_{\cal P}$ induces an inner product on  $T_u{\cal A}^L_{\cal P}$ as subspace of $T_u{\cal C}^L_{\cal P}$.   In fact these inner products are proportional and moreover,  the norm $||\;||_\infty$ and $||\;||_{L^2}$ induce equivalent norm on  $T_u{\cal A}^L_{\cal P}$.
\end{rem}

\subsection{The horizontal distribution associated to a Hilbert snake}
For any $u\in {\cal C}^L_{\cal P}$ consider the map $S_u:[0,L]\ap \field{H}$ given by:
\begin{eqnarray}\label{serpent}
S_u(t)=\dis\int_0^tu(s)ds
\end{eqnarray}
called the {\bf Hilbert snake} associated to $u$.  On the other hand, to each configuration $u\in {\cal C}^L_{\cal P}$ we can associate the {\bf endpoint map}:
\begin{eqnarray}\label{end}
{\cal E}: & {\cal C}^L_{\cal P} \ap \field{H}\nonumber\\
& \qquad  u\ap S_u(L)
\end{eqnarray}

As ${\cal E}$ is the restriction to ${\cal C}^L_{\cal P}$ of the linear map $u\ap \dis\int_0^Lu(s)ds$ defined on ${\cal C}^0_{\cal P}([0,L],\field{H})$ it follows that ${\cal E}$ is smooth and we have:
\begin{eqnarray}\label{tangeps}
T_u{\cal E}(v)=\dis\int_0^Lv(s)ds
\end{eqnarray}
Note that   ${\cal E}[{\cal C}^L_{\cal P}]$  is the closed ball $\field{B}_L=\{x\in \field{H} \textrm{ such that } ||x||\leq L\}$.

\begin{lem}\label{horozontal} ${}$
\begin{enumerate}
\item The subspace $\ker T_u{\cal E}\subset T_u{\cal C}^L_{\cal P}$  is a Banach space for each induced norm $||.||_\infty$ and $||.||_{L^2}$.
\item  The orthogonal of $\ker T_u{\cal E}$ (for the inner product $<.,.>_{L^2}$ on $T_u{\cal C}^L_{\cal P}$), denoted   ${\cal D}_u$, is a closed  space in each normed spaces  $ (T_u{\cal C}^L_{\cal P},||.||_{L^2})$ and $ (T_u{\cal C}^L_{\cal P},||.||_\infty)$ and  we have the decomposition
\begin{eqnarray}\label{horidecomp}
 T_u{\cal C}^L_{\cal P}={\cal D}_u\oplus \ker T_u{\cal E}.
 \end{eqnarray}
 \item In the Banach space  $ (T_u{\cal C}^L_{\cal P},||.||_\infty)$, the restriction of $T_u{\cal E}$ to ${\cal D}_u$ is a continuous  injective morphism into $\field{H}$
\end{enumerate}
\end{lem}

 \textit{Proof }:\\
 At first,  it is well known  that we have (see\cite{Die} (8.7.7)) that:
$$||T_u{\cal E}(v)||\leq\dis\int_0^L||v(s)||ds$$
So on one hand we get
$$||T_u{\cal E}(v)||\leq L||v||_\infty$$
On the other hand by Cauchy-Schwartz inequality we get
$$||T_u{\cal E}(v)||\leq \sqrt{L}||v||_{L^2}$$

\noindent So $\ker T_u{\cal E}$ is closed in $T_u{\cal C}^L_{\cal P}$ according to each  norm $||.||_\infty$. and $||.||_{L^2}$, which ends the proof of part {\it 1.}\\

By construction,  ${\cal D}_u$ is a closed subset of  the normed space $ (T_u{\cal C}^L_{\cal P},||.||_{L^2})$; so we have the decomposition  (\ref{horidecomp}) in this norm space. On the other hand, it follows from (\ref{inegalitenorm}) that ${\cal D}_u$  is also closed in the Banach space  $ (T_u{\cal C}^L_{\cal P},||.||_\infty)$ and  then, the decomposition  (\ref{horidecomp}) is again true in this Banach space which ends the proof of part {\it 2.}\\

\noindent According  to the decomposition (\ref{horidecomp}) in the Banach space  $ (T_u{\cal C}^L_{\cal P},||.||_\infty)$ we get part {\it 3.} \\
${}$\hfill$\D$\\

\begin{defi} \label{horizon} ${}$
\begin{enumerate}
\item The family $u \mapsto {\cal D} _u$ is a (closed) distribution on ${\cal C}^L_{\cal P}$ called  the {\bf  horizontal distribution}.
\item Each vector field $X$ on ${\cal C}^L_{\cal P}$ which is tangent to ${\cal D}$ is called a {\bf horizontal vector field}.
\end{enumerate}
\end{defi}

On ${\cal A}^L_{\cal P}$, the intersection ${\cal D}_u\cap T_u{\cal A}^L_{\cal P}$ gives rise to a (closed) Hilbert distribution ${\cal D}^{\cal A}$. Note that we can also define ${\cal D}^{\cal A}$ directly as the orthogonal  of $\ker T_u{\cal E}\cap T_u{\cal A}^L_{\cal P}$ relatively to one of the equivalent inner products defined on $T_u{\cal A}^L_{\cal P}$ (see Remark \ref{weakmetric}). When no confusion is possible, this distribution ${\cal D}^{\cal A}$  on ${\cal A}^L_{\cal P}$ will be also denoted by $\cal D$ and also called the {\bf horizontal distribution on  ${\cal A}^L_{\cal P}$}.\\

The inner product on $\field{H}$ gives rise to a Riemannian  metric $g$ on $T\field{H}\equiv \field{H}\times\field{H}$ given by $g_x(u,v)=<u,v>$.
Let $\phi:\field{H}\ap \R$ be a smooth function. The usual gradient  of $\phi$ on $\field{H}$ is the vector field
$$\mathrm{grad}(\phi)=(g^\f)^{-1}(d\phi)$$
where $g^\f$ is the canonical isomorphism of bundle  from $T\field{H}$ to its dual bundle $T^*\field{H}$, corresponding to the Riesz representation i.e.
$g^\f(v)(w)=<v,w>$. So $\mathrm{grad}(\phi)$ is characterized by:
\begin{eqnarray}\label{grad}
 g( \mathrm{grad}(\phi),v)=<\mathrm{grad}(\phi),v>=d\phi(v)
\end{eqnarray}
for any $v\in \field{H}$.

On the opposite, on  $T{\cal C}^L_{\cal P}$,  the Riemannian metric ${G}$ is {\it only  weak} (see Remark \ref{weakmetric}) and we cannot define in the same way the gradient of any smooth function on  ${\cal C}^L_{\cal P}$. However let be \\${G}^\f:T{\cal C}^L_{\cal P}\ap T^*{\cal C}^L_{\cal P}$ the morphism bundle  defined by:
$${G}_u^\f(v)(w)={G}_u(v,w)$$
for any $v$ and $w$ in $T_u{\cal C}^L_{\cal P}$. Then we have

\begin{lem}\label{gradsurC}${}$
 Let $\phi: \field{H}\ap \R$ be a smooth function. $\ker d(\phi\circ {\cal E})$  contains $\ker T{\cal E}$ and belongs to $\tilde{G}_u^\f(T_u{\cal C}^L_{\cal P})$. Moreover,
\begin{eqnarray}
\nabla\phi=({G}^\f)^{-1}(d(\phi\circ{\cal E}))
\end{eqnarray}
 is tangent to ${\cal D}_u$. Moreover, we have
\begin{eqnarray}\label{expnabla}
\nabla\phi(u)(s)=\mathrm{grad}(\phi)({\cal E}(u))-<\mathrm{grad}(\phi)({\cal E}(u)),u(s)>u(s)
\end{eqnarray}
\end{lem}

\begin{rem}\label{finitgrad}${}$
When $\field{H}$ is finite dimensional, the relation (\ref{expnabla}) is exactly  the definition of $\nabla\phi$ given in \cite{Ro}.
\end{rem}
\bigskip
\begin{defi}${}$\\
For any smooth function $\phi: \field{H}\ap \R$, the vector field  $\nabla\phi$ is called  {\bf horizontal gradient} of $\phi$.
\end{defi}
\textit{Proof Lemma \ref{gradsurC}}: given $v\in T_u{\cal C}^L_{\cal P}$ then we have:
$d(\phi\circ{\cal E})(v)=d\phi(T{\cal E}(v))=d\phi(\dis\int_0^Lv(s)ds)$.
\noindent So if $v\in \ker T_u{\cal E}$ then $d(\phi\circ{\cal E})(v)=0$. So, according to the decomposition (\ref{horidecomp}), it follows that  we get
\begin{eqnarray}\label{reduc}
d(\phi\circ{\cal E})(v)=d(\phi\circ{\cal E})(\pi_u(v))
\end{eqnarray}
where $\pi_u:T_u{\cal C}^L_{\cal P}\ap{\cal D}_u$ is the canonical projection associated to the decomposition  (\ref{horidecomp}).

On the other hand, consider the vector field  $$\mathrm{Grad}(\phi)(u)(s)=\mathrm{grad}(\phi)({\cal E}(u))-<\mathrm{grad}(\phi)({\cal E}(u)),u(s)>u(s)$$ along  $u$.  As we have

$$<\mathrm{Grad}(\phi)(u)(s),u(s)>=0$$ it follows that $\mathrm{Grad}(\phi)(u)$ belongs to $T_u{\cal C}^L_{\cal P}$. Moreover, for any $v\in T_u{\cal C}^L_{\cal P}$, as $<u(s),v(s)>=0$ for any $s\in [0,L]$, we get

${G}(\mathrm{Grad}(\phi)(u),v)=\dis\int_0^L<\mathrm{grad}(\phi)({\cal E}(u)),v(s)>ds$

${}\qquad\qquad \qquad\quad\;\;=<\mathrm{grad}(\phi)({\cal E}(u)),\dis\int_0^Lv(s)ds>$

${}\qquad\qquad \qquad\quad\;\;=d\phi({\cal E}(u))[\dis\int_0^Lv(s)ds]$

${}\qquad\qquad \qquad\quad\;\;=d\phi({\cal E}(u))\circ T_u{\cal E}(v)$

The first consequence of the last equality is that $\mathrm{Grad}(\phi)(u)$ belongs to ${\cal D}_u$. Using the identity $d\phi({\cal E}(u))\circ T_u{\cal E}(v)=d(\phi\circ {\cal E})(u)(v)$, the second consequence is that $d(\phi\circ {\cal E})$ belongs ${G}_u^\f(T_u{\cal C}^L_{\cal P})$ which ends the proof of Lemma \ref{gradsurC}

${}$\hfill $\D$

To each vector $x\in \field{H}$, we can associate the linear form $x^*$ such that $x^*(z)=<z,x>$. So from Lemma \ref{gradsurC} the horizontal gradient  $\nabla x^*$ is well defined. In particular, to each vector  $e_i$, $i\in \N$, of the Hilbert basis, we can associate the horizontal vector field $E_i=\nabla e_i^*$. Then  as in \cite{Ro} we have:

\begin{lem}\label{baseD}${}$\\
The family $\{E_i\}_{i\in\N}$ of vector fields generates the distribution ${\cal D}$.
\end{lem}

\textit{ Proof}:   Let $u\in {\cal C}^L_{\cal P}$ be; we can write
$$u(s)=\dis\sum_{i\in \N}u_i(s)e_i$$
  Denote by $\D_u$ the closed subspace generated by the family $\{E_i(u)\}_{i\in \N}$ in the normed space $(T_u{\cal C}^L_{\cal P},||.||_{L^2})$. A vector $v\in T_u{\cal C}^L_{\cal P}$ belongs to the orthogonal of $\D_u$ (relatively to $G$) if and only if $G(v,E_i(u))=0$ for all $i\in N$. But as $<v(s),u(s)>=0$ we have:
\begin{eqnarray}\label{GH}
G(v,E_i)=\dis\int_0^L<v(s),e_i- u_i(s)u(s)>ds=<\dis\int_0^L v(s),e_i>\textrm { for all } i\in \N
\end{eqnarray}
According to (\ref{GH})  $v$ is orthogonal to $\D_u$ if and only if  $v\in \ker T_u{\cal E}$. As ${\cal D}_u$  is also closed in $(T_u{\cal C}^L_{\cal P},||.||_{L^2})$ , we get $\D_u={\cal D}_u$.\\
${}$\hfill$\D$

\begin{rem}\label{surS}${}$
\begin{enumerate}
\item As in finite dimension (see \cite{Ro}), for $\phi=e^*_i$  using the left member of (\ref{expnabla}), for any $i\in \N$  we have

$E_i(s)=e_i-<e_i,u(s)>u(s)$

So,  each $E_i(u)$ can be considered as a vector field on $\field{S}^\infty$ along $u:[0,L]\ap \field{S}^\infty$. In this way, $E_i(u)$ is nothing but the orthogonal projection of $e_i$ onto the tangent space to $S^\infty$ along $u([0,L])$.
\item On ${\cal A}^L_{\cal P}$ the induced inner product $<\;,\;>_{L^2}$ induces a (strong) Riemannian metric on the horizontal distribution $\cal D$. In the same way, to the Hilbert basis $\{e_i,\;i\in \N\}$ of $\field{H}$ we can associate a family of global vector fields (again denoted) $\{E_i,\;i\in \N\}$ on ${\cal A}^L_{\cal P}$. In fact these vector fields are only the restriction to ${\cal A}^L_{\cal P}$ of the family defined on the whole manifold ${\cal C}^L_{\cal P}$.
if we identify $T{\cal A}^L_{\cal P}$ with $ [T\field{S}^\infty]^N$ (see Remark \ref{weakmetric}), the vector field $E_i$ at $u=(x_1,\cdots,x_N)$ is
$$( e_i-<x_1,e_i>x_1,\cdots e_i-<x_N,e_1>x_N).$$
 If  there is no ambiguity, we also  denote these family in the same way. Of course, on ${\cal A}^L_{\cal P}$, the distribution $\cal D$ is also generated by this family of vector fields.
\end{enumerate}
\end{rem}

\subsection{Set of critical values and set of singular points of the endpoint map}\label{sing}

As the continuous linear map $T_u{\cal E}: T_u{\cal C}^L_{\cal P}\ap T_{{\cal E}(u)}\field{H}\equiv \field{H}$ is closed (see \cite{Die} section 8.7), it follows that $\r_u=T_u{\cal E}_{|{\cal D}_u}$ is  an isomorphism from ${\cal D}_u$ to  the closed subset $\r_u({\cal D}_u)$ of $\field{H}$. Consider a point  $u\in {\cal C}^L_{\cal P}$. According to remarks \ref{weakmetric} {\it 1.}, the annulator of $\r_u({\cal D}_u)=T_u{\cal E}({\cal D}_u)$ is
 $$[\r_u({\cal D}_u)]^0=\{z\in T_{{\cal E}(u)}\field{H}\equiv \field{H}\textrm{ such that }<z, \r_u(v)>=0, \;\forall v\in {\cal D}_u\}$$
So, $u$ is a {\bf  singular point} of $\cal E$ if and only if
$[\r_u({\cal D}_u)]^0\not=\{0\}$ .\\

\noindent On the other hand, as the family $\{E_i(u)\}_{i\in \N}$ generates ${\cal D}_u$, any $z\in \field {H}$ belongs to  $[\r_u({\cal D}_u)]^0$ if and only if we have
\begin{eqnarray}\label{sing1}
<z, \dis\int_{0}^{L} E_i(u)(s)ds>=0,\; \forall i\in \N
\end{eqnarray}

\noindent Consider the decompositions $z=\dis\sum_{i\in \N}z_ie_i$ and $u(s)=\dis\sum_{i\in \N}u_i(s)e_i$. Then (\ref{sing1}) is equivalent to
\begin{eqnarray}\label{sing2}
Lz_i=\dis\sum_{j\in \N}\int_0^Lu_i(s)u_j(s)z_jds\; \forall i\in \N
\end{eqnarray}
Let $\G_u$ be the endomorphism  defined by matrix of general term $ (\int_0^Lu_i(s)u_j(s)ds)$. Note that $\G_u$ is self-adjoint. The endomorphism $A_u=L.Id-\G_u$ is also self-adjoint and, in fact,
 its matrix  in the basis $\{e_i\}_{i\in \N}$ is $ (L\d_{ij}-\int_0^Lu_i(s)u_j(s)ds)$. So (\ref{sing2}) is equivalent to
\begin{eqnarray}\label{sing3}
A_u(z)=0
\end{eqnarray}
 So $u$ is a singular point if and only if $L$ is an eigenvalue of $\G_u$.  The proof of the following Lemma \ref{4.4} is an adaptation of the argument used  in finite dimension  (see\cite{Ro})

\begin{prop}\label{pointsing}${}$\\
A point  $u\in {\cal C}^L_{\cal P}$ is a singular point of $\cal E$ in and only if  the vector space generated by $u([0,L])$ is $1$-dimensional.
\end{prop}
\smallskip

\textit{ Summarized proof}: at first, note that for any unitary automorphism $U$ of $\field{H}$ we have
\begin{eqnarray}\label{invUG}
U\G_uU^*=\G_{U(u)}
\end{eqnarray}
On the other hand, we have
\begin{eqnarray}\label{invUE}
{\cal E}\circ U(u)=U({\cal E}(u))
\end{eqnarray}
 If $u([0,L])$ generates a $1$-dimensional space then we have $u(s)=\pm x\in \field{S}^\infty$. Using (\ref{invUG}), without loss of generality, we can suppose that $u(s)=\pm e_1$ for any $s\in[0,L]$. In this case, using the  relation obtained  by derivation of  (\ref{invUE}) we show that  $e_1$ is an eigenvector associated to the eigenvalue $L$ of $\G_u$ and so $\ker (L.Id-\G_u)=\ker A_u\not=\{0\}$.

 On the other hand, if $u$ is a singular point of $\cal E$, there exists a vector $x\in \field{S}^\infty$ which is an eigenvector associated to the eigenvalue $L$ of $\G_u$. If $U$ is an unitary automorphism such that $U(x)=e_1$ then $e_1$ is an eigenvector associated to $L$ for $U\G_uU^*=\G_{U(u)}$. If we set $\bar{u}=U(u)$ then we get
 $\G_{\bar{u}}(e_1)=Le_1$. So, for the decomposition $\bar{u}(s)=\dis\sum_{i\in \N}\bar{u}_i(s)e_i$, we get  $\dis\int_0^L[\bar{u}_1]^2=L$ and $\bar{u}_i(s)\equiv 0$ for all $i>1$. It follows that  $\bar{u}(s)=\pm e_1$ and so $u(s)=\pm x$.\\
 ${}$\hfill$\D$\\

According to Lemma  \ref{pointsing},  a point $u\in{\cal C}^L_{\cal P}$  is  singular if and only if  the restriction to $[s_{i-1},s_i]$ is equal to $\pm x$ for some $x\in \field{S}^\infty$.
  It follows that  the {\bf set of singular points} $\S({\cal E})$ of ${\cal E}$ is diffeomorphic   to the projective space  a  $\field{P}^\infty$ of $\field{H}$. \\

 On the other hand, let $u\in\S({\cal E})$ be with $u(s)\equiv x\in \field{S}^\infty$. For any $v\in {\cal C}^0_{\cal P}([0,L],\field{H}) \textrm { such that }\\ <u(s),v(s)>=0 \textrm{ for all } s\in [0,L]$ we consider
 $$\bar{u}_n=\dis\frac{1}{n}v(s)+x$$
 As $||x||=1$, for $n$ large enough, we have $||\bar{u}_n(s)||\geq \dis\frac{1}{2}$ and  $u_n(s)=\dis\frac{\bar{u}_n}{||\bar{u}_n(s)||}$ belongs to ${\cal C}^L_{\cal P}\setminus \S({\cal E})$. Moreover we have
 $$\dis\lim_{n\ap\infty}u_n=u$$
 \noindent So the set ${\cal C}^L_{\cal P}\setminus \S({\cal E})$ of  {\bf regular points } of $\cal E$  is an open  dense subset of ${\cal C}^L_{\cal P}$.\\

Recall that the image of $\cal E$ is the closed ball $B(0,L)$ in $\field{H}$. As in finite dimension,  when ${\cal P}=\{0,L\}$ the set of critical values of $\cal E$ is then  the  boundary of $B(0,L)$ i.e. the  sphere $S(0,L)$ and $\{0\}$. In the general case, ${\cal P}= \left\{ a=s_{0}<s_{1}<...<s_{N}=b\right\} $, the same argument applied to each subinterval $[s_{i-1},s_i]$ gives that the {\bf set of critical values} of $\cal E$ is the union of spheres $S(0,L_j)$ for $j=1,\cdots n$  with $0\leq L_j\leq L$.

\begin{rem}\label{hilbertstruc}${}$
\begin{enumerate}
\item Recall that $ \r_u$ is an isomorphism from ${\cal D}_u$ to  $\r_u({\cal D}_u)$, which is a closed subspace of $\field{H}$.  So  on ${\cal D}_u$, the  norm induced by $||.||_\infty$ is equivalent  the norm  $||.||_{\field{H}}$; moreover,  $\r_u$ is an isometry between ${\cal D}_u$ and $\r_u({\cal D}_u)$ endowed with the Hilbert induced norm. In particular, for  any  regular  point $u$, the inverse of $\r_u$ is given by $\dis\frac{1}{L}\nabla v^*$
and according to (\ref{GH}) we have $\r_u(\dis\frac{1}{L}\nabla v^*)=v$. So $\{E_i(u),\;i\in \N\}$ is then a Hilbert basis  of ${\cal D}_u$  according to this isometry.
If now, $u$ is a singular point of $\cal E$, according to the previous proof, there exists a Hilbert basis $\{e'_i,\;i\in \N\}$ of $\field{H}$ such that $e'_1=u(s)$ for all $s\in [0,L]$. So, $\r_u$ is an isomorphism from ${\cal D}_u$ to $\{e'_1\}^\perp$. It follows that, on ${\cal D}_u$, the  norm induced by $||.||_\infty$ is equivalent to the norm $||.||_{\field{H}}$ so that $\r_u$ is an isometry  between ${\cal D}_u$ and  $\{e'_1\}^\perp$. Then,  the family $\{E'_i(u)=\nabla (e'_i)^*(u),\; i>1\}$ is a Hilbert basis of ${\cal D}_u$.
\item  According to the beginning of this section, as $G(E_i(u), E_j(u))=L\d_{ij}-\int_0^Lu_i(s)u_j(s)ds$, the matrix of $G$ in the basis $\{E_i(u),\; i\in \N\}$ is the matrix of $A_u=L.Id- \Gamma_u$. But $A_u$ is a self-adjoint endomorphism of $\field{H}$ which is compact. So the sequence $\{\l_i,\;i\in \N\}$ of eigenvalues of $A_u$ is bounded and converges to $0$ and there exists a Hilbert basis $\{e'_i,\;i\in \N\}$ of eigenvectors of $A_u$. In this basis, the matrix of $A_u$ is diagonal and equal to  $(L-\l_i\d_{ij})$.  So for the associated family $\{E'_i(u),\;i\in \N\}$ of generators of ${\cal D}_u$ we have:

(1) if $u$ is regular the matrix of $G$  in the basis $\{E'_i(u),\;i\in \N\}$ is $(L-\l_i)\d_{ij}$. Note that $0$ is not an eigenvalue of $A_u$ otherwise, it would  mean that  $u$ is an eigenvector of $\G_u$ associated to the eigenvalue $L$. As  the sequence  $\{\l_i,\;i\in \N\}$ is bounded and converges to $0$, there exists $K>0$ so that   $\frac{1}{K}\leq L-\l_i\leq K$  for any $i\in \N$. It follows that    the norm associated to $G$ and the norm associated to the isometry $\r_u$ are equivalents.

(2) if $u$ is singular, according to the proof of Lemma  \ref{pointsing}, we can choose $e'_1$ so that $u=\pm e'_1$  and then, by the same arguments as the ones used in (1) but  applied to the restriction of $A_u$ to $\{e'_1\}^\perp$ we again  obtain that the norm associated to $G$ and the norm associated to the isometry $\r_u$ are equivalent.
\end{enumerate}
\end{rem}

Finally we obtain the following result :

\begin{prop}\label{sub}${}$
\begin{enumerate}
\item The set  ${\cal R}({\cal E})$ (resp. ${\cal V({\cal E})}$) of regular values (resp. points) of $\cal E$ is an open dense subset of ${\cal C}^L_{\cal P}$ (resp. $\field{H}$).
\item
 For any $u\in{\cal R}({\cal E})$  the linear map $\r_u:{\cal D}_u\ap \{{\cal E}(u)\}\times \field{H}$   is an isomorphism and  on ${\cal D}_u$,  the inner product induced by $<\;,\;>_{L^2}$) and the inner product
defined  $\r_u$ from $\field{H}$ are equivalent.  Moreover the distribution  ${\cal D}_{|{\cal R}({\cal E})}$ defines trivial  Banach bundle over ${\cal R}({\cal E})$
\item  The distribution  ${\cal D}_{|{\S({\cal E})}}$ is an Hilbert  bundle which  is isometrically  isomorphic  to $T\field{P}^\infty$
\end{enumerate}
\end{prop}
\noindent\so{\it Proof}: we consider the map $F: {\cal C}^L_{\cal P}\times l^2(\N)\ap {\cal C}^L_{\cal P}\times {\cal C}^0_{\cal P}([0,L],\field{H})$ defined by
$$F(u,\s)=\dis\sum_{i\in \N}\s_iE_i(u)$$
It is easy to see that
 the family of smooth vector fields $\{E_i,\; i\in \N\}$ satisfies the condition (LBs) for any $s\in \N$ at any point  and as  $F$ is linear in $\s$, it follows that $F$ is a smooth map. According to  Lemma \ref{baseD} and Remark \ref{hilbertstruc}, the range of $F$ is ${\cal D}$.  Again, from Remark \ref{hilbertstruc}, $\{E_i(u),\; i\in \N\}$ is an Hilbert basis of ${\cal D}_u$ for  $u\in {\cal R}({\cal E})$. So, the restriction of $F$ to ${\cal R}({\cal E})\times l^2(\N)$ is a global trivialization of  ${\cal D}_{|{\cal R}({\cal E})}$. The same argument can be used for the restriction of $F$ to $\S({\cal E})$. The other claims were proved previously.\\
${}$\hfill $\D$

\section{ The optimal control for a Hilbert snake}
\subsection{The problem of optimality and accessibility}\label{result}

As we have seen in the introduction, recall that  our  {\bf optimality problem} for a Hilbert snake can be formulated in the following way:
 given any $C^k$-piecewise map $\g: [0,T] \ap\field{H}$ we look for a $1$-parameter  family $\{u_t\}_{t\in [0,1]}$  such that the associated  family $S_t$ of snakes satisfies  $S_t(L)=\g(t)$ for all $t\in [0,1]$ so that the infinitesimal kinematic energy $\dis\frac{1}{2}\int_0^L||\frac{d}{dt}u_t(s)||ds$ is minimal. \\

 Given the horizontal distribution $\cal D$ on ${\cal C}^L_{\cal P}$, a {\bf horizontal curve} is $C^k$-piecewise map \\$\g: [0,T] \ap {\cal C}^L_{\cal P}$ such that $\dot{\g}(s)$ belongs to ${\cal D}_{\g(s)}$ almost everywhere. On the other hand, given any $C^k$-piecewise curve $c:[0,T]\ap \field{H}$, a {\bf lift} of $c$ is a $C^k$-piecewise curve $\g: [0,T] \ap {\cal C}^L_{\cal P}$ such that
${\cal E}(\g(t))=c(t)$. When a lift $\g$ is horizontal we say that $\g$ is a {\bf horizontal lift}. By construction of ${\cal D}$, among all lifts of $c$, a horizontal lift is a lift which minimizes infinitesimally the kinematic energy $\dis\frac{1}{2}G(\dot{\g},\dot{\g})$. So our optimal problem has a solution if and only if the curve $c$ has a horizontal lift. On the other hand, we can also ask when two positions $x_0$ and $x_1$  of the "head" of the snake can be joined  by a piecewise smooth curve $c$ which has an "optimal control" ${\g}$ as lift. As in finite dimension, the {\bf accessibility set} ${\cal A}(u)$, for some $u\in{\cal C}^L_{\cal P}$,  is the set of endpoints $\g(T)$ for any piecewise smooth   horizontal curve $\g:[0,T]\ap{\cal C}^L_{\cal P} $ such that $\g(0)=u$. In this case if $x_0=S_u(L)$ then any $z=S_{u'}(L)$ can be joined from $x_0$ by an absolutely continuous  curve   $c$ which has an "optimal control" when $u'$ belongs to ${\cal A}(u)$.\\


In finite dimension, given any horizontal distribution ${\cal D}$ on a finite dimension manifold $M$, the famous {\it Sussmann's Theorem} (see \cite{Su}) asserts that  each accessibility set is a smooth immersed manifold which is an integral manifold of a distribution $\hat{\cal D}$ which contains $\cal D$  (i.e. ${\cal D}_x\subset \hat{\cal D}_x$ for any $x\in M$) and  characterized by:

 $\hat{\cal D}$ is the smallest distribution which contains $\cal D$  and which is invariant by the flow of any (local) vector field tangent to $\cal D$.\\

In the context of Banach manifolds  the reader can find some generalization of this result in \cite{LaPe}. In the next section, we will use the results of this paper to give some  positive answer to this accessibility problem via an  analogue construction as previously.  More precisely, according to subsection  \ref{orbit}, we will  associate to each Hilbert basis $\{e_i,\;i\in\N\}$ of $\field{H}$,  the family ${\cal X}=\{E_i,\;i\in \N\}$  of (global) vector fields (see Lemma \ref{baseD}) which can be extended to a family $\{E_i, [E_j,E_k],\;i,j,k\in \N, k<l\}$ such that the associated distribution $\hat{\cal D}$ is integrable. Moreover, this last set of vector fields generates a weak distribution $\bar{\cal D}$ modeled on Hilbert spaces  with the following properties:

(i) $\bar{\cal D}$ does not depend on the choice of the basis  $\{e_i,\;i\in\N\}$;

(ii) $\hat{\cal D}_x$ is dense in $\bar{\cal D}_x$ for all $x\in M$;

(iii) $\bar{\cal D}$ is integrable and each maximal integral manifold of $\bar{\cal D}$ contains the orbit of $\{E_i,\;i\in \N\}$ for any choice of basis $\{e_i,\;i\in\N\}$ of $\field{H}$;

(iv) the accessibility set of any point of such a maximal integral $N$  manifold is a dense subset of  $N$.\\

In this way we  obtain:
\begin{theo}\label{acesset}${}$
Let  $\{e_i,\; i\in \N\}$ be a Hilbert basis of $\field{H}$ and $\{E_i,\; i\in \N\}$ the associate family of vector fields on  ${\cal C}^L_{\cal P}$. The vector space
$$\bar{\cal D}_u=\{\dis\sum_{i\in \N}x_iE_i(u)+\sum_{j,l\in \N,j<l} \xi_{ij}[E_i,E_j](u)\,, \sum(x_i)^2<\infty,\sum (\xi_{ij})^2<\infty\}$$
  is a well defined subspace of $T_u{\cal C}^L_{\cal P}$ and carries a natural structure of Hilbert space such that the inclusion of $\bar{\cal D}_u$ in $T_u{\cal C}^L_{\cal P}$ is continuous and gives rise to a weak Hilbert  distribution on ${\cal C}^L_{\cal P}$. This distribution has the following properties:
\begin{enumerate}
\item[(1)]  $\bar{\cal D}$ does not depend on the choice of the Hilbert basis $\{e_i\}$ of $\field{H}$.
\item[(2)] The distribution  $\bar{\cal D}$  is integrable. Moreover,  for each  $u\in {\cal C}^L_{\cal P} $, the accessibility set ${\cal A}(u)$ is a dense subset of the maximal  integral manifold $L(u)$ of $\hat{\cal D}$  through $u$
\item[(3)] on the manifold ${\cal A}^L_{P}$, each subspace  $\bar{\cal D}_u\cap T_u{\cal A}^L_{P}$ induces a closed distribution (again denoted by $\bar{\cal D}$) which satisfies the two previous properties and  moreover, in this case, each integral maximal manifold of this distribution is a Hilbert submanifold of ${\cal A}^L_{P}$.
\end{enumerate}
 \end{theo}
\begin{rem}\label{min}${}$\\
Recall that a horizontal curve $\g$ is an absolutely continuous curve in ${\cal C}^L_{\cal P}$ which is almost everywhere  tangent to $\cal D$.
Given $u\in {\cal C}^L_{\cal P}$, we denote by $H_u\subset T_u {\cal C}^L_{\cal P}$ the set of tangent vectors at $u$ of a horizontal curve through $u$ which has a tangent vector at $u$. If $X$ and $Y$ are vector fields on $ {\cal C}^L_{\cal P}$ whose domain contains $u$, the curve
$$t\ap \Phi^X_t\circ\Phi^Y_t\circ\Phi^X_{-t}\circ\Phi^Y_{-t}(u)$$
is a horizontal curve and it is well known that its tangent vector at $u$ is $[X,Y](u)$. So, if we look for the smallest (weak) manifold of ${\cal C}^L_{\cal P}$ which contains the accessibility set ${\cal A}(u)$, its tangent space must contain $H_u$. In particular, this tangent space must contain the family $\{E_i(u),[E_j,E_l](u), i,j,l\in \N\}$. Note that from Theorem \ref{acesset}, it follows that $\bar{\cal D}_u$ contains $H_u$. On one hand, if we consider the closed distribution generated by ${\cal X}=\{E_i,[E_j,E_l], i,j,l\in \N\}$, we can show that this distribution is upper trivial and the property (1) of  Theorem \ref{acesset} is satisfied. But we do not know if this distribution is integrable. On the other hand, according to the following subsection, the $l^1$- weak distribution $\D^1$ generated by ${\cal X}$ satisfies property (2), but not property (1) and so $\D^1_u$ does not contain $H_u$. So, in this sense the distribution $\bar{\cal D}$ is the "smallest"  weak distribution which is integrable and such that the  maximal integral manifold through $u$ contains ${\cal A}(u)$. Moreover, as, the maximal integral manifold $N$ is closed in this case, the ${\cal X}$-orbit of $u$ is contained in $N$,  for any family ${\cal X}$ of type $\{E'_i,\; i\in \N\}$ associated to any Hilbert basis $\{e'_i,\;i\in \N\}$ of $\field{H}$.
On the other hand,  when $\field{ H}$ is finite dimensional, $\bar{\cal D}$ is exactly the distribution whose leaves are  the accessibility sets as proved in \cite{Ro}.
\end{rem}

According to  our problem of optimality for the head of the snake, we know that if $u$ is a configuration,  and $N$ is the maximal integral manifold  of $\bar{\cal D}$ through $u$, for all other  configuration $v\in N$ there exists a sequences $(\g_n)$ of horizontal curves in $N$ whose  origin $u$ and whose sequence extremities converges to $v$. So, if ${\cal E}(u)=x$ and ${\cal E}(v)=y$, the family of curves $c_n={\cal E}\circ \g_n$ are optimal (in the previous sense),  have $x$ for origin, and,  the sequence of extremities $y_n$ of $c_n$ converges to $y$.

For each  maximal integral manifold $N$ of $\bar{\cal D}$, denote by $\tilde{N}$ the  range $\tilde{N}={\cal E}(N)$. Then for each pair  $(x,y)\in \tilde{N}$ there exists a family of optimal  curves $c_n$  which  have $x$ for origin, and,  the sequence of extremities $y_n$ of $c_n$ converges to $y$. \\

\subsection{Construction of the distribution $\bar{\cal D}$}\label{barD}

For the construction of $\bar{\cal D}$ we need the following result whose proof  is the same as in the case of a finite dimensional Hilbert space $\field{H}$ (see \cite{Ro}). According to Remark \ref{surS}, each $E_i$ can be considered as a vector field on $\field{S}^\infty$. In these way, we have

\begin{lem}\label{relbrack}${}$\\
The brackets of vector fields of the family $\{E_{i}\}_{i\in \N}$ satisfy the following relations:

$[E_i,E_j](u)=<e_j,u>E_i(u)-<e_i,u>E_j(u)$ for any $u\in{\cal C}^L_{\cal P}$ and any $i,j\in \N$;

$[E_i[E_j,E_k]]=\d_{ij}E_k-\d_{ik}Ej$  for any $i,j,k\in \N$

 $[[E_i,E_j],[E_k,E_l]]=\d_{il}[E_j,E_k]+\d_{jk}[E_i,E_l]-\d_{ik}[E_j,E_l]-\d_{jl}[E_i,E_k]$ for any $i,j,k,l\in \N$.\\
\end{lem}

\noindent We consider  the countable  set of indexes $\L=\{(i,j),\; i,j\in \N,\; i<j\}$ and let $\field{G}^1$ (resp $\field{G}^2$)  be  the Banach space $l^1(\N)\oplus l^1(\L)$ (resp. $l^2(\N)\oplus l^2(\L)$). We then have the following result:


\begin{lem}\label{4.4}${}$
 \begin{enumerate}
  \item For $p=1,2$,   the map $\Psi^p$ from the trivial bundle ${\cal C}^L_{\cal P}\times \field{G}^p$ to $T{\cal C}^L_{\cal P}$ characterized by:
\begin{eqnarray}\label{defPsi}
\Psi_u^p(\s,\xi)=\dis\sum_{i\in \N}\s_iE_i(u)+ \sum_{(i,j)\in \L}\xi_{ij}[E_i,E_j](u),\; \s=( \s_i)\in l^p(\N),\;  \xi=(\xi_{ij})\in l^p(\L)
\end{eqnarray}
is well defined and each $\Psi_u^p$ is a  continuous linear map.
\item  For each $u\in {\cal C}^L_{\cal P}$, let  $\field{V}_u$ be the Hilbert subspace of $\field{H}$ generated by the set
$$\{u(t)-u(0),\; t\in [0,L]\}$$
 For $p=1,2$,  if the kernel of $\Psi_u^p$ is not $\{0\}$ then $\field{V}_u\not=\field{H}$
\item For $p=1,2$, the distribution $\D^p$ defined by ${\D}^p_u=\Psi^p_u(\field{G}^p)$ is a weak distribution and the map  $\Psi^p$ defines a strong (global) upper  trivialization of ${\D}^p$.
\item The distribution ${\D}^2$  do not depend of the choice the Hilbert basis $(e_i)$ in $\field{H}$  contains ${\cal D}$ (i.e. ${\cal D}_u\subset{\D}^2_u$)
\end{enumerate}
\end{lem}
\bigskip

\noindent\so{ \it Proof of  Lemma \ref{4.4}}\\

\noindent\so{\it Proof of  part 1}${}$\\
For any $\s\in l^p(\N)$  the vector $\dis\sum_{i\in \N}\s_ie_i$ belongs to $\field{H}$ and for any $s\in[0,L]$ the vector $\dis\sum_{i\in \N}\s_iE_i(u(s))$ is the orthogonal projection on $T_{u(s)}\field{S}^\infty$ of  $\dis\sum_{i\in \N}\s_ie_i$. So we have
\begin{eqnarray}\label{proj}
||\dis\sum_{i\in \N}\s_iE_i(u)||_\infty\leq[\dis\sum_{i\in \N}(\s_i)^2]^{1/2}=||\s||_{2}
\end{eqnarray}
If $\s$ belongs to $l^1(\N)$, as $||\s||_{2}\leq |\s||_1$ in this case we get
$$||\dis\sum_{i\in \N}\s_iE_i(u)||_\infty\leq||\s||_1$$
On the other hand as  $[E_k,E_l](u)=u_lE_k(u)-u_kE_l(u)$, in the same way, for any $s\in [0,L]$, the vector $\sum_{(k,l)\in \L}\xi_{kl}[E_k,E_l](u(s)) $ is the orthogonal projection of $\dis \sum_{(k,l)\in \L}\xi_{kl}(u_l(s)e_k-u_k(s)e_l)$ on  $T_{u(s)}\field{S}^\infty$.

\noindent But we have:
$$\dis \sum_{(k,l)\in \L}\xi_{kl}u_l(s)e_k=\dis\sum_{k\in \N}[\sum_{l>k}\xi_{kl}u_l(s)]e_k]$$
so
$$||\dis \sum_{(k,l)\in \L}\xi_{kl}u_l(s)e_k||^2=\dis\sum_{k\in \N}[\sum_{l>k}\xi_{kl}u_l](s)]^2$$
Using the fact that $|u_j(s)|\leq ||u(s)||=1$,  from Cauchy -Schwartz inequality we get:
$$|\dis\sum_{l>k}\xi_{kl}u_l](s)|\leq [\sum_{l>k}(\xi_{kl})^2]^{1/2}$$
Finally we obtain
$$||\dis \sum_{(k,l)\in \L}\xi_{kl}u_l(s)e_k||^2\leq \dis \sum_{(k,l)\in \L}(\xi_{kl})^2=(||\xi||_2)^2$$
By same argument we get
$$||\dis \sum_{(k,l)\in \L}\xi_{kl}u_k(s)e_l ||^2\leq \dis \sum_{(k,l)\in \L}(\xi_{kj})^2=(||\xi||_2)^2$$
So we obtain
$$||\sum_{(k,l)\in \L}\xi_{kl}[E_k,E_l](u)||_\infty\leq 2||\xi||_2 $$

If $\xi\in l^1(\L)$  by same argument as previously we also get:
$$||\sum_{(k,l)\in \L}\xi_{kl}[E_k,E_l](u)||_\infty\leq 2||\xi||_1 $$
 Finally  we get :
 \begin{eqnarray}\label{cont}
||\Psi^p_u(\s,\xi)||_\infty\leq 2||(\s,\xi)||_p \textrm{ for } p=1,2
\end{eqnarray}

It follows that $\Psi^p$ is well defined. From its expression, it is easy to see that $\Psi^p_u$ is linear and continuous from (\ref{cont}). This ends the proof of part 1.\\

\noindent\so{\it Proof of part 2}${}$\\
At first, note that as the natural  inclusion $I:\field{G}^1\hookrightarrow \field{G}^2$ is continuous and with dense range, we have $\Psi^2_u\circ I=\Psi^1_u$,  the closure of  $\ker \Psi^1_u$ in $\field{G}^2$ is equal to $\ker \Psi^2_u$. So $\ker \Psi^1_u\not=0$ if and only if   $\ker \Psi^2_u\not=0$.  Assume that $\ker \Psi^1_u\not=\{0\}$\\

Let be $(\s,\xi)\in\ker \Psi^1_u$. According to (\ref{defPsi}), and Remark \ref{surS}, we must have

\begin{eqnarray}\label{ker1}
\dis\sum_{i\in \N}\s_iE_i(u(s)) + \sum_{(i,j)\in \L}\xi_{ij} [u_i(s)E_j(u(s))-u_j(s)E_i(s)]=0 \textrm{ for any } s\in [0,L]
\end{eqnarray}
  We set  $\bar{\xi}_{kj}=\dis\frac{\xi_{kj}}{2}$ (resp. $\bar{\xi}_{kj}=-\dis\frac{\xi_{kj}}{2}$) for $j<k$ (resp $j>k$) and $\bar{\xi}_{jj}=0$. Then
(\ref{ker1}) can be written:

\begin{eqnarray}\label{ker2}
\dis\sum_{i\in \N}[\sum_{j\in \N}(\bar{\xi}_{ij}u_j(s)+\s_i]E_i(u(s))=0 \textrm{ for any } s\in [0,L]
\end{eqnarray}
Given any $\xi\in l^1(\L)$, denote  $\Xi$ the endomorphism of $l^1(\N)$ whose matrix in the canonical basis is precisely $(\bar{\xi}_{ij},\;{i,j\in \N})$. So, (\ref{ker2}) is equivalent to
\begin{eqnarray}\label{ker3}
\Xi u(s)=-\s  \textrm{ for any } s\in [0,L]
\end{eqnarray}
So, $\s$ must belong to the range of $\Xi$.



\noindent According to the definition of $\field{V}_u$,  (\ref{ker3}) is equivalent to

 $\Xi u(0)=-\s$ and  $\field{V}_u\subset \ker \Xi$.\\

\noindent\so{\it Proof of part 3}${}$\\
From part 1,  $\D^p_u=\Phi^p_u(\field{G}^p)$ gives rise to a well defined distribution on ${\cal C}^L_{\cal P}$. On the other hand, denote by $\hat{\Psi}^p_u$ the canonical bijection induced by $\Psi^p_u$:
$$\hat{\Psi}^p_u:\field{G}^p/\ker \Psi^p_u\ap \D^p_u$$
So we can put on $\D^p_u$ the Banach structure so that $\hat{\Psi}^p_u$ is an isometry. In this way, $\D^p$ is then a weak distribution. On the other hand, the family
of  smooth vector field $\{E_i, \;i\in \N\}$ satisfies the condition (LBs) for any $s\in \N$ at any point, and as $\Psi^1_u$ is linear with Lipschitz constant independent of $u$,  the map $(u,(\s,\xi)) \mapsto \Psi^p_u(\s,\xi)$ is smooth.\\

It remains to show that $\ker \Psi^p_u$ is complemented in $\field{G}^p$ for each $u\in {\cal C}^L_{\cal P}$. At first, for $p=2$, as $\field{G}^2$ is a Hilbert space, it is always true. In particular, the previous Banach structure on each $\D^2_u$ is a Hilbert structure. However, we shall show this result for each case  $p=1$ and $p=2$.\\

{\it assume that $u\in {\cal R}({\cal E})$}\\
If $\ker\Psi^p_u=\{0\}$ there is nothing to prove. Now assume that  $\ker\Psi^p_u\not=\{0\}$.
 At first, suppose that  we have a partition $\N=A\cup B$ such that $\{e_a,\; a\in A\}$ (resp. $\{e_b,\;b\in B\})$ is a Hilbert basis of $[\field{V}_u]^\perp$ (resp. $\field{V}_u)$.  By construction,  each component $u_a$ is constant, for all  $a\in A$.  
  So the Lie brackets $[E_a,E_{a'}]$, for $a,a'\in A$, belongs to ${\cal D}_u$. Let be
$$\field{K}=\{\xi\in l^p(\L) \textrm{ such that } \xi_{ij}=0 \textrm{ if } i \textrm{ or } j \in B\}$$
According to the notations of the proof of part 2, for any $\xi\in \field{K}$ if we denote again by $\Xi$ the  associated endomorphism of $\field{H}$, $\ker \Xi$ contains $\field{V}_u$ and  if $\s=-\Xi u(0)$ then $(\s,\xi)$ belongs to $\ker \Psi^p_u$. So, the subspace
$$\hat{\field{K}}=\{(\s,\xi)\in l^p(\N)\oplus l^p(\L),\; \xi\in \field{K}, \s=-\Xi u(0)\}$$
is contained in $\ker\Psi^p_u$.

\noindent On the other hand, if $(\s,\xi)$ belongs to $\ker \Psi^p_u$, from the proof of part 2  and we have $\s=-\Xi u(0)$  and $\field{V}_u\subset \ker \Xi$ and, as $\{e_b\}$ is a basis of $\field{V}_u$, we then have    $(\s,\xi)\in \hat{\field{K}}$. It follows that $\ker \Psi^1_u$ is complemented:

 if we denote by   $\field{L}$ is the subspace of
$\{\xi\in l^1(\L),\; \xi_{ij}=0 \textrm{ for all } i,j\in A\}$,
then the subspace $l^p(\N)\oplus \field{L}$ is a complement subspace of $\ker\Psi^p_u$.\\

In the general case, choose a Hilbert basis $\{e'_a,\; a\in A\}$ (resp $e'_b,\; b\in B\}$) of $[\field{V}_u]^\perp$ (resp. $\field{V}_u$).
There exists a linear isometry $T$ of $\field{H}$ such that $T(e'_a)=e_a$ for $a\in A$ and $T(e'_b)=e_b$ for $b\in B$. Denote by $E'_j=\nabla (e'_j)^*$ the associated vector field on ${\cal C}^L_{\cal P}$ (see Lemma \ref{gradsurC}). The map $\tilde{T}:(z, v)\ap (z,T(v))$ is an isomorphism of $T\field{H}$ such that $\tilde{T}(E'_j)(u)=E_j(u)$ for any $j\in \N$.
Consider the map  $\Psi': \field{G}^1\times {\cal C}^L_{\cal P}\ap T{\cal C}^L_{\cal P}$ defined by
$$\Psi'_u(\s,\xi)=\dis\sum_{i\in \N}\s_iE'_i(u)+ \sum_{(i,j)\in \L}\xi_{ij}[E'_i,E'_j](u),\; \s=( \s_i)\in l^p(\N),\;  \xi=(\xi_{ij})\in l^p(\L)$$
Of course, we have
$$\Psi^p_u=\Psi'_{u}\circ T$$
But in the new basis, for $\Psi'_u$ we are in the previous situation.
So, it follows that $\ker \Psi^p_u$ is complemented, which ends the proof of part 3. \\

{\it Assume now that $u\in \S({\cal E})$}.\\
 According to the proof of Lemma \ref{pointsing}, then $u(t)=\pm x\in\field{S}^\infty$ for any $t\in [0,L]$, there exists an Hilbert basis $\{e'_i,\;i\in \N\}$ such that $x= e'_1$,  the associated family $\{E'_i(u),\; i>1\}$ is  a basis of ${\D}_u$ and we have $E'_1(u)=0$ (see Remark \ref{baseD}). Moreover, as the components of $u$ are constant, from Lemma \ref{relbrack}, all brackets $[E'_j,E'_l](u)$ belongs to ${\cal D}_u$ for $i>1$ and $j>1$ and $[E'_1,E'j]=-x_iE_j$ also belongs to  ${\D}_u$. As previously, we can consider the map
 $$\Psi'_u(\s,\xi)=\dis\sum_{i\in \N}\s_iE'_i(u)+ \sum_{(k,l)\in \L}\xi_{kj}[E'_i,E'_j](u),\; \s=( \s_i)\in l^p(\N),\;  \xi=(\xi_{ij})\in l^p(\L)$$
  Its kernel is $\R e'_1\oplus l^p(\L)$. From the same argument as previously, we obtain that the $\ker \Psi^p_u$ is complemented. Moreover, the restriction of $\Psi^p_u$ to $l^P(\N)$ has a kernel of dimension $1$ and the restriction of $\Psi^p$ to the orthogonal $[\ker \Psi^p_u]$ in $l^1(\N)$ is an isomorphism onto  $\D^p_u$.\\

\noindent\so{\it Proof of part 4}${}$\\
 Now, we must  show that the range of $\Psi^2_u$ does not depend on the choice of the Hilbert basis $(e_i)$ of $\field{H}$. So given any other Hilbert basis $(e'_j)$ of $\field{H}$, denote again  by  $E'_j=\nabla (e'_j)^*$ the associated vector field on ${\cal C}^L_{\cal P}$. So we have the decomposition:
 \begin{eqnarray}\label{decompo}
E'_i=\dis\sum_{j\in \N}a_i^j E_j \textrm{ and }
[E'_j,E'_k]=\dis\sum_{l,m \in \N}a_j^l a_k^m[E_l,E_m]=\sum_{(l,m)\in\L}(a_j^la_k^m- a_j^ma_k^l)[E_l,E_m]
\end{eqnarray}

Let $T$ be the isometry of $\field{H}$ defined by $T(e_i)=e'_i$. Let $l^2_{BA}(\field{H})$ be the set of Hilbert-Schmidt  bilinear antisymmetric maps. Then $\{e^*_i\wedge e^*_j,\; (i,j)\in \L\}$ is a Hilbert basis  of $l_{BA}^2(\field{H})$. Denote by $T^2$ the isometry  of $l_{BA}^2(\field{H})$ induced by $T$ on $l^2_{BA}(\field{H})$. Then matrix of $T^2$ in this basis is precisely $[(a_i^ka_j^l- a_j^ka_i^l)]_{(i,j), (k,l)\in \L}$. It follows that $T$ (resp. $T^2$) is an isometry of $l^2(\N)$ (resp. $l^2(\L))$

 On the other hand, to the choice $(e'_i)$ of a basis of $\field{H}$ is naturally associated   the map ${\Psi'}^2:{\cal C}^L_{\cal P}\times\field{G}^p\ap T{\cal C}^L_{\cal P}$ characterized by:
   $${\Psi'}^2_u(\s,\xi)=\dis\sum_{i\in \N}\s_iE'_a(u)+ \sum_{(j,k)\in \L}\xi_{jk}[E'_j,E'_k](u)$$
  According to (\ref{decompo}) we have:
  $$\Psi^2_u(\s,\xi)={\Psi'}^2_u(T\s,T^2\xi)$$
   So we have ${\Psi'}^2_u(\field{G}^2)=\D^2_u$  for any $u\in{\cal C}^L_{\cal P}$. \\

 On the other hand, according to  Lemma \ref{baseD}, it is clear that ${\cal D}_u$ is contained in ${\D}^2_u$.  According to Remark \ref{hilbertstruc}, We can note that $\psi^2_u(l^2(\N)={\cal D}_u$  and from the proof of part 3, $\Psi^2(l^2(\L))$ is a complemented space of ${\cal D}_u$ in $\D^2_u$\\

${}$\hfill$\D$

\subsection{Proof of Theorem \ref{acesset}}\label{access}
According to Lemma \ref{4.4} , we have  $\bar{\cal D}_u=\D^2_u$ and so $\bar{\cal D}$ is a well defined weak Hilbert distribution on ${\cal C}^L_{\cal P}$ which does not depend on the choice of the Hilbert basis $\{e_i,\; i\in \N\}$ of $\field{H}$.\\
We take place in the context of the proof of Lemma \ref{4.4}. According to Lemma \ref{relbrack}, Lemma \ref{4.4} and Theorem \ref{lieinv}, it follows that the distribution $\D^1$ is integrable. On the other hand, according to Lemma \ref{relbrack}, Lemma \ref{4.4}  and Theorem \ref{lieinv} we also have that
$\bar{\cal D}=\D^2$ is also integrable.
We again denote by $\cal X$ the family $\{E_i,\;i\in \N\}$ of vector fields.\\
{\it Now We will show that any $\cal X$-orbit is contained in a maximal integral manifold of $\bar{\cal D}$.}\\

For the sake of simplicity, we only denote by $\field{G}$ the previous Hilbert space $\field{G}^2$.
As the distribution $\bar{\cal D}$ is integrable, let $f:N\ap {\cal C}^L_{\cal P}$ be any maximal integral manifold of $\bar{\cal D}$. Without loss of generality, we can identify $N$ with $f(N)$ and take $f=i_N$ the natural inclusion of $N$ (with its Hilbert manifold structure) into ${\cal C}^L_{\cal P}$.   Consider the pull-back $f^*( {\cal C}^L_{\cal P}\times  \field{G}) $ over $N$. Note that $f^*({\cal C}^L_{\cal P}\times  \field{G})$ can be identified with $N\times \field{G}$. As the range of $\Psi_u$ is $\bar{\cal D}_u$, for any $u$, the bundle morphism $\Psi: {\cal C}^L_{\cal P}\times  \field{G}\ap T{\cal C}^L_{\cal P}$ induces a bundle morphism $\tilde{\Psi}$ from  $ N\times \field{G}$ to $TN$ which is onto.  Moreover,  the orthogonal of $\ker \tilde{\Psi}_u$ in $\{u\}\times \field{G}$, gives rises to a Hilbert sub-bundle of $N\times \field{G}$. Denote by $\cal N$ this sub-bundle and by ${\Pi}$ the natural orthogonal projection of $N\times \field{G}$ on $\cal N$. Now, we have $\Pi\circ\tilde{\Psi}=\tilde{\Psi}$ and the restriction of $\tilde{\Psi}$ to $\cal N$ is an isomorphism  from $\cal N$ onto $TN$ and we have
\begin{eqnarray}\label{compoPsi}
Tf\circ\tilde{\Psi}= \Psi\circ(Id\times f)
\end{eqnarray}
 Now, given the canonical  Hilbert basis $\{\e_i,\o_{jl},\; i\in \N, (j,l)\in \L\}$ of $\field{G}$, we set $\hat{E}_i(u)=\tilde{\Psi}_u(\e_i)$ and $\hat{E}_{jl}(u)=\tilde{\Psi}_u(\o_{jl})$. In fact, $\hat{E}_i$ and $\hat{E}_{jl}$ are smooth global vector fields on $N$.
On the other hand,  we have $\Psi_v(\e_i)(u)=E_i(u)$ and $\Psi_v(\o_{jl})=[E_j,E_l](v)$ for any $v\in f(N)$.\\  According to proposition \ref{varXS}, there exist (global)  vector fields $\tilde{E}_i$ on $N$ such that $f_*\tilde{E}_i=E_i$ and so $f_*[\tilde{E}_j,\tilde{E}_l]=[E_j,E_l]$. It follows from (\ref{compoPsi})  that
$\hat{E}_i=\tilde{E}_i$ and $\hat{E}_{jl}=[\tilde{E}_j,\tilde{E}_l]$.

Let $\tilde{\cal X}$ be the induced family $\{\tilde{E}_i,[\tilde{E}_j,\tilde{E}_l],\; i,j,l\in \N,\; j<l\}$.
 As $\Psi$ (resp. $\tilde{\Psi}$) is a strong (global) upper trivialization for $\bar{\cal D}$ on ${\cal C}^L _{\cal P}$ (resp $TN$ on $N$), it follows that, for any $u\in {\cal C}^L _{\cal P}$ (resp.$u\in N$), there exists an open neighborhood $U\subset  {\cal C}^L _{\cal P}$ (resp.  $\tilde{U}\subset N $) of $u$ such that $\cal X$  (resp. $\tilde{\cal X}$) satisfies the condition (LBs) on $U$  (resp. $\tilde{U}$) for $s>3$ (see \cite{LaPe}, proof of Theorem 6, part 2).\\
Consider any family $\xi=\{X_\a,\;\a\in A\}\subset{\cal X}$ and let $\tilde{\xi}=\{\tilde{X}_\a,\; \a\in A\}$ be the corresponding family on  a maximal integral manifold $N$.  Given  $u\in f(N)$, consider some flow  $\Phi^\xi_\t$ associated to $\xi$ and let $\g(t)=\Phi^\xi_\t(t,u)$ be the integral curve defined on $[0,||\t||_1]$. From proposition  \ref{varXS}, there exists a curve $\tilde{\g}: [0,||\t||_1[\ap N$ such that $f\circ \tilde{\g}=\g$ on $[0,||\t||_1[$. We set $v=\g(||\t||_1)$.\\

{\it We want to show that $v$ also belongs to $f(N)$, or equivalently, $\Phi^\xi_\t(||\t||_1,u)=\phi^\xi_\t(u)$ belongs to $N$.}\\

 Consider a maximal integral manifold  $g\equiv i_M:M\ap {\cal C}^L_{\cal P}$ of $\bar{\cal D}$ through $v$ and set $\tilde{v}=(i_M)^{-1}(v)$. As, we have already seen, if ${\cal X}'$ is the family of vector fields $\{E'_i,[E'_i,E'_l],\;i,j,l,\in \N, j<l\}$ on $M$ such that $g_*E'_i=E_i$, then  ${\cal X}'$ satisfies the condition (LBs). On $N$, we also have a family  ${\xi'}=\{X'_\a,\a\in A\}$ defined on a neighborhood of $v$ and so
$g_*X'_\a=X_\a$. Then $\xi'$ also satisfies the condition (LBs) for $s>3$. So, from Theorem 2 of \cite{LaPe}, there exists $\eta>0$ such that,  for $\t'\in l^1(A)$  with $||\t'||_1\leq \eta$, the corresponding flow $\Phi^{\xi'}_{\t'}(.,.)$ is defined on a neighborhood $\tilde{V}$ of $\tilde{v}=g^{-1}(v)$ in $M$. Now coming back to the original flow $\Phi^\xi_\t$ on ${\cal C}^L_{\cal P}$, if $\t=(\t_\a)_{\a\in A}$, there exists $\a_0$ such that

$\dis\sum_{\a\geq \a_0}|\t_\a|< \eta$

Then, given any $a\in A$ with $a\geq \a_0$, we set   $\t_a=(\t'_\a)$ with $\t'_\a=0$ for $\a<a$ and $\t'_\a=\t_\a$ for $\a\geq a$. The corresponding flows $\Phi^{\xi'}_{\t_a}$ and $\hat{\Phi}^{\xi'}_{\t_a}$ are defined on $M$. Moreover, we have
\begin{eqnarray}\label{xi'}
g\circ \Phi^{\xi'}_{\t_a}(t,\tilde{z})=\Phi^\xi_{\t_a}(t,g(\tilde{z}))\textrm{ and } g\circ\hat{ \Phi}^{\xi'}_{\t_a}(t,\tilde{z})=\hat{\Phi}^\xi_{\t_a}(t,g(\tilde{z}))
\end{eqnarray}
for any $\tilde{z}\in \tilde{V}$.\\
By construction of the flow $\Phi^\xi_\t$ we have

$\Phi^\xi_{\t_a}(||\t_a||_1,\g(\t_a))=v$ and so $\hat{\Phi}^\xi_{\t_a}(||\t_a||_1, v)=\g(\t_a)$

For any  $a\geq \a_0$, in ${\cal C}^L_{\cal P}$ consider the curve  $\hat{\g}_a(s)=\hat{\Phi}^\xi_{t_a}(||\t_a||_1-s, v)$. This curve is defined on $[0,||\t_a||_1]$ and  joins $v$ to $\g(\t_a)$. In the  same way, in $M$, consider the curve $\hat{\g}'_a(s)=\hat{\Phi}^{\xi'}_{\t_a}(||\t_a||_1-s, v)$. This curve is also defined on  $[0,||\t_a||_1]$  and joins $\tilde{v}$ to $\tilde{v}_a$ in $N$. According to (\ref{xi'}) we have
$$g\circ\hat{\g}'_a=\hat{\g}_a.$$
 In particular, we get  $g(\tilde{v})=\g(\t_a)$. But $\g(\t_a)$ belongs to $f(N)\equiv N$  and to $g(M)\equiv M$ as subsets of  ${\cal C}^L_{\cal P}$. But, $(N,f\equiv i_N)$ and $(M,g\equiv i_M)$ are maximal integral manifolds of $\bar{\cal D}$. So, as $N\cap M\not=\emptyset$,  we must have  $N=M$ and so we can extend $\tilde{\g}$ to the closed interval $[0,||\t||_1]$ and, in particular, $\phi^\xi_\t(u)=\Phi^\xi_\t(||\t||_1,u)$ belongs to $N$.\\

Now if  we have $v=\Phi(u)$ for some $\Phi\in {\cal G}_{\cal X}$ (see subsection \ref{orbit}), then $\Phi$ is a  finite composition of local diffeomorphisms  of type $\phi^\xi_\t$ or $[\phi^\xi_\t]^{-1}$ or of  type $\Phi^X_t$ for some $X\in {\cal X}$. From the previous argument, if $u\in L$, then  $\phi^\xi_\t(u)$ and $[\phi^\xi_\t]^{-1}(u)$ belong to $N$ and from Lemma \ref{varXS}, part 1,  $\Phi^X_t(u)$ also belongs to $N$. By induction we obtain that $v=\Phi(u)$ belongs to $N$. So the $\cal X$-orbit ${\cal O}(u)$ of $u$ is contained in $N$.\\

{\it If ${\cal A}(u)$ denotes the accessibility set of  $u$ we now show that ${\cal A}(u)$ is dense in $N$.}\\

From  Theorem \ref{III"},  ${\cal O}(u)$  contains the maximal integral manifold $L^1$ of $\D^1$ through $u$ and $L^1$ is dense in ${\cal O}(u)$ (for the topology of ${\cal C}^L_{\cal P}$, and so  $L^1$ and ${\cal O}(u)$ have the same closure in ${\cal C}^L_{\cal P}$.  But  we have $L^1\subset {\cal O}(u)\subset L$.  As the inclusion of $N$ in  ${\cal C}^L_{\cal P}$ is continuous we obtain that ${\cal O}(u)$ is dense in  $N$.

From proposition  \ref{l1Xorbit}  the set ${\cal A}(u)\cap {\cal O}(u)$ is dense in  ${\cal O}(u)$. On the other hand, according to proposition  \ref{varXS}, we see that ${\cal A}(u)$ is contained in $N$.  So we obtain that ${\cal A}(u)$ is dense in $N$.\\

Now, it remains to show that the same results are true on ${\cal A}^L_{\cal P}$.  It is easy to see (and it is left to the reader)  that all the proofs of Lemma \ref{4.4}  work in the same way on the manifold ${\cal A}^L_{\cal P}$. So the previous arguments work too in this context. But, in ${\cal A}^L_{\cal P}$, the corresponding  distribution $\bar{\cal D}$ is closed. So each maximal integral manifold of $\bar{\cal D}$ is a weak Hilbert manifold whose topology is the topology induced by the topology of the Hilbert manifold ${\cal A}^L_{\cal P}$. So, such a manifold must be a Hilbert submanifold of ${\cal A}^L_{\cal P}$.

${}$\hfill $\D$
\subsection{Almost Lie  algebroid  structures }\label{ALalgebroid}
According to Lemma \ref{4.4},
for $p=1,2$, on $\field{G}^p$ we define a  Lie algebra structure in the following way:

let be $(\epsilon_i)_{i\in \N}$ (resp. $(\epsilon_{ij})_{(i,j)\in \L}$ the canonical basis of $(l^p(\N)$ (resp. $(l^p(\L))$;

according to Lemma \ref{relbrack}, we then define:

$[\e_i,\e_j]=\o_{ij}$, for all $i,j\in \N$

$[\e_i,\o_{jk}]=\d_{ij}\e_k-\d_{ik}\e_j$, for all $i\in \N$ and $(j,k)\in \L$

 $[\o_{ij},\o_{kl}]=\d_{il}\o_{jk}+\d_{jk}\o_{il}-\d_{ik}\o_{jl}-\d_{jl}\o_{ik}$ , for all $(i,j) (kl)\in \L$.\\

\noindent For any $x=\dis\sum x_i\a_i$, $y=\sum y_j\e_j$ in $l^p(\N)$ and  $\xi=\sum \xi_{ij}\o_{ij}$,  $\eta=\sum\eta_{kl}\b_{kl}$ in $l^p(\L)$, naturally we can define:

$[\s,\s']=\dis\sum_{i,j\in \N}\s_i\s'_j[\e_i,\e_j]$

$[\s,\eta]=\dis\sum_{i\in \N, (k,l)\in \L}\s_i\eta_{kl}[\e_i,\o_{kl}]$

 $[\xi,\eta]=\dis\sum_{(i,j)\in \L,(k,l)\in \L}\xi_{ij}\eta_{kl}[\o_{ij},\o_{kl}]$.\\

 Now, according to Lemma \ref{4.4},  the map $\Psi^p:{\cal C}^L_{\cal P}\times \field{G}^p\ap T{\cal C}^L_{\cal P}$ is morphism bundle over ${\cal C}^L_{\cal P}$. Moreover, each section $\varphi$  of the trivial bundle ${\cal C}^L_{\cal P}\times \field{G}\ap {\cal C}^L_{\cal P}$ can be identified with a map $\varphi:{\cal C}^L_{\cal P}\ap \field{G}^p$. So, on the set $\G(\field{G}^p)$ of section of this trivial bundle we can defined a Lie bracket by:

$$[\varphi,\varphi'](u)=[\varphi(u),\varphi'(u)]+d\varphi(\Psi^p(u,\varphi'(u))-d\varphi'(\Psi^p{u,\varphi(u)})$$
 According to \cite{Pe} section 4 or \cite{CaPe},  it follows that  $(\field{G}^p\times {\cal C}^L_{\cal P}\,\Psi,{\cal C}^L_{\cal P},[\;,\;])$ has  a {\bf  Banach  Lie algebroid structure}  on ${\cal C}^L_{\cal P}$\\

In $\field{G}^p$ let  be $\pi:\field{G}^p\ap l^p(\N)$ the canonical projection  whose kernel is $l^p(\L)$ and denote again by $\pi:{\cal C}^L_{\cal P}\times  \field{G}^p\ap {\cal C}^L_{\cal P}\times l^p(\N)$ the associated  projection bundle. Again  any section of the trivial bundle ${\cal C}^L_{\cal P}\times l^p(\N)\ap{\cal C}^L_{\cal P} $ can be identified with a map from ${\cal C}^L_{\cal P}$ to $l^p(\N)$. Of course the set $\G(l^p(\N))$  of such sections is contained in $\G(\field{G}^p)$. So, according to  \cite{CaPe}, on $\G(l^p(\N))$, we can define an almost Banach  Lie bracket by:
$$[[\varphi,\varphi']](u)=\pi([\varphi,\varphi'](u)).$$
So, if we denote by  $\theta^p$ the restriction of $\Psi^p$ to ${\cal C}^L_{\cal P}\times l^p(\N)$
we get an {\bf almost  Banach Lie algebroid structure} $({\cal C}^L_{\cal P}\times l^p(\N),\Psi^,{\cal C}^L_{\cal P},[\;,\;])$ on ${\cal C}^L_{\cal P}$.

\noindent Moreover, again, according to \cite{CaPe} subsection 4.3,    the inner product on $l^2(\N)$ gives rise to  {\bf strong Riemaniann metric} on  $({\cal C}^L_{\cal P}\times l^2(\N),\Psi^2,{\cal C}^L_{\cal P},[\;,\;])$. Note that, according to Remark \ref{hilbertstruc}, the induced inner product on ${\cal D}_u$ is equivalent to the inner product associated to the Riemaniann metric $G$.

  Given a maximal   integral manifold   $(f,N)$ of $\bar{\cal D}$,  the pull back $f_*( {\cal C}^L_{\cal P}\times l^2(\N))$ and $f^*( {\cal C}^L_{\cal P}\times \field{G}^2)$ can be identified with $N\times l^2(\N)$ and $ N\times \field{G}$ respectively. Then,  $\theta^2$ and $\Psi^2$ induces  anchors   $\theta_N:N\times l^2(\N)\ap TN$  and $\Psi_N: N\times \field{G}^2\ap TN$ and the almost bracket   $[[\;,\;]]$ induces an almost bracket again denoted $[[\;,\;]]$. So $(N\times  l^2(\N),\theta_N,N,[[\;,\;]])$ is an {\it  almost  Banach Lie algebroid} on $N$  and     $(N\times  \field{G}^2,\Psi_N,N,[[\;,\;]])$ is a  {\it  Banach Lie algebroid} on $N$.  Moreover the canonical scalar product on $l^2(\N)$ (resp. $\field{G}^2$) gives rise to a strong Riemannian metric on  $(N\times  l^2(\N),\theta_N,N,[[\;,\;]])$ (resp. on $(N\times  \field{G}^2,\Psi_N,N,[[\;,\;]])$).

 Essentially from the proof of Lemma \ref{4.4}  we get:

   \begin{prop} \label{PropN}${}$\\
 Fix some  $u\in N$. Then we have the following properties
  \begin{enumerate}
  \item[(1)] Assume that  $u\in \S({\cal E})$ then $N=\S({\cal E})$. Let be ${\cal L}_v$  the $1$-codimensional Hilbert  subspace  $[\ker \theta_N]_v^\perp \subset \{v\}\times l^2(\N)$ for any $v\in N$. Then ${\cal L}=\dis\cup _{v\in N}{\cal L}_v$ is a $1$-codimensional Hilbert  sub-bundle of $N\times l^2(\N)$ and the restriction $\psi_N$ of $ \Psi_N$ to ${\cal L}$ is an isomorphism onto $TN$ and we have ${\cal D}_{| N}=TN$.
   \item [(2)] Assume that $u\in {\cal R}({\cal E})$.
   Let be  $\field{V}_u$ the Hilbert subspace of $\field{H}$ generated by the set
$$\{u(t)-u(0),\; t\in [0,L]\}$$
 and choose an Hilbert   basis $\{e'_a,\; a\in A\}$ (resp $e'_b,\; b\in B\}$) of $[\field{V}_u]^\perp$ (resp. $\field{V}_u$) (see the proof of part 3 of Lemma \ref{4.4} ). If $\L_u$ is the set of pair $(i,j)\in \L$ such that
 that $i$ or $j$ do not belongs to $A$, then $N$ is an Hilbert manifold modeled on $l^2(\N)\oplus l^2(\L_u)$  and is contained in $ {\cal R}({\cal E})$ .\\
  Let be ${\cal L}_v$ the orthogonal of $\ker [\Psi_N]_v\subset   \{v\}\times \field{G}^2$. Then ${\cal L}=\dis\cup _{v\in N}{\cal L}_v$ is a Hilbert  sub-bundle of $N\times \field{G}^2$   which contains $N\times l^2(\N)$ and the restriction of $\psi_N$ of $\Psi_N$ to ${\cal L}$ is an isomorphism on $TN$ \\
  Moreover, ${\cal L}$ contains $N\times l^2(N)$ and the restriction  of $\theta_N$ to $N\times l^2(N)$ is an isomorphism on ${\cal D}_{| N}$
 \item[(3)] Let be $\g$  an horizontal  piecewise $C^k$ curve in $N$. Then there exists an unique  piecewise $C^{k-1}$ section  $(\g,\s)$ of $\cal L$ over $\g$  such that $\psi_N(\g(t)),\s(t))=\dot{\g}(t)$ for all $t$.
 \end{enumerate}
 \end{prop}

 \noindent\so{\it Proof}

\noindent  From the proof of part 3 of Lemma \ref{4.4}  for $p=2$  , we get  that $\bar{\cal D}_u={\cal D}_u $ for any  $u\in \S({\cal E})$. So $\S({\cal E})$ is an integral manifold of $\bar{D}$.
 On the other hand,  from the proof of part 3 of Lemma \ref{4.4}    for $p=2$, taking  $u\in {\cal R}({\cal E})$ we get that the tangent space to $N$ is modeled on  $l^2(\N)\oplus l^2(\L_u)$  (with the  notations introduced in  part  (2)). As for  $u\in {\cal R}({\cal E})$, we have ${\cal D}\not=\bar{\cal D}$, it follows that  $\S({\cal E})$ is a maximal integral manifold of $\bar{D}$, and, for $u\in {\cal R}({\cal E})$, so  $N$ is contained in ${\cal R}({\cal E})$.\\

As $\Psi_N:N\times \field{G}^2\ap TN$ is a surjective Hilbert  bundle morphism, the kernel of this morphism is an Hilbert  sub-bundle ${\cal K}$ of $N\times  \field{G}^2$. It follows that ${\cal L}$ is also an Hilbert sub-bundle of $N\times \field{G}^2$ and  the restriction $\psi_N$ of $\Psi_N$ to each ${\cal L}$ is an isomorphism onto $TN$.  Given $u\in N$, if $u\in {\cal R}({\cal E})$, according to notations in (2),  by the same arguments used in proof of part 4 of Lemma \ref{4.4}, we can show that  $l^2(\N)$ is contained in ${\cal L}_u$. Now, if we identify $\field{H}$ with $l^2(\N)$, for  $u\in\S({\cal E})$,  the kernel of $\Psi^2_u$ in $\{u\}\times \field{G}^2$ is $\R.{\cal E}(u)\oplus l^2(\L)$. So the vector space ${\cal L}_u$ is the orthogonal of $\ker [\Psi_N]_u$ in $l^2(\N)$ (see the proof of part 3 of Lemma \ref{4.4} ).\\

Let be $\g:[0,T]\ap N$ a piecewise $C^k$ horizontal curve  for $k \geq 1$. Assume that $N\subset{\cal R}({\cal E})$.
 On one hand, $\psi_N:{\cal L}\ap TN$ is an isomorphism and on the other hand $\psi_N(N\times l^2(\N))={\cal D}_{| N}$ for  $u\in {\cal R}({\cal E})$ (resp. $\psi_N({\cal L}= {\cal D}_{| N}$ for $u\in \S({\cal E})$. So the restriction $\theta_N$ of  $\psi_N$ to $N\times l^2(\N)$  is an isomorphism of bundle onto ${\cal D}_{| N}$. As $\g$ is horizontal, the curve  $(\g(t),x(t))=(\g(t),[\theta_N]^{-1}(\dot{g}(t))$ is a well defined piecewise $C^{k-1}$ curve which satisfies the conclusion in (3). Finally, when  $N=\S({\cal E})$,  $\psi_N$ is an isomorphism from $\cal L$ on $TN$ (see part (1))so conclusion (3) is clear in this case.

 ${}$\hfill $\D$

\end{document}